\documentclass{article}

\usepackage[utf8]{inputenc}
\usepackage{enumitem}
\usepackage{setspace}
\usepackage{mathtools}
\usepackage{tikz-cd}
\usepackage{lipsum}
\usepackage{tikz}
\usepackage{mathrsfs}
\usetikzlibrary{arrows}
\usepackage{amssymb}
\usetikzlibrary{decorations.pathmorphing}
\usepackage{bbold}
\usepackage{amsmath,amsthm}
\usepackage{xspace}
\usepackage{faktor}
\usepackage{calrsfs}
\usepackage{verbatim}
\usepackage{quiver}
\usepackage[pdfa]{hyperref}

\newtheorem{theorem}{Theorem}[section]
\newtheorem{proposition}[theorem]{Proposition}
\newtheorem{lemma}[theorem]{Lemma}

\newtheorem{corollary}[theorem]{Corollary}

\newtheorem{definition}[theorem]{Definition}
\newtheorem{notation}[theorem]{Notation}

\newtheorem{remark}[theorem]{Remark}

\newcommand{\lgrp}{\ell \mathbb{G}rp}

\newcommand{\cat}[1]{\mathbb{#1}}
\newcommand{\mono}{\rightarrowtail}
\newcommand{\inv}[1]{#1^{-1}}
\newcommand{\pol}[1]{#1^{\perp}}
\newcommand{\splt}[2]{\mathbb{S}pl \mathbb{E}xt_{#2}(#1)}
\newcommand{\cev}[1]{\reflectbox{\ensuremath{\vec{\reflectbox{\ensuremath{#1}}}}}}

\newcommand{\conj}[2]{#1^{-1}#2#1}
\newcommand{\jnoc}[2]{#1#2#1^{-1}}

\newcommand{\Pt}[2]{Pt_{#1}\mathbb{#2}}
\newcommand{\pt}[2]{Pt_{#1}{#2}}
\newcommand{\grp}{\mathbb{G}rp}
\newcommand{\lab}{\ell \mathbb{A}b}

\DeclareMathOperator{\kerr}{ker}
\DeclareMathOperator{\Aut}{Aut}
\DeclareMathOperator{\Imm}{Im}
\DeclareMathOperator{\point}{Pt}
\DeclareMathOperator{\Eq}{Eq}
\DeclareMathOperator{\Act}{SplExt}
\DeclareMathOperator{\coeq}{coeq}
\DeclareMathOperator{\Ideals}{Ideals}

\DeclareMathOperator{\terminal}{\boldsymbol{1}}

\usepackage[english]{babel}	

\title{Categorical-Algebraic Properties of Lattice-ordered Groups}
\author{Andrea Cappelletti}

\begin{document}

\maketitle

\begin{abstract}

We study the categorical-algebraic properties of the semi-abelian variety $\lgrp$ of lattice-ordered groups. In particular, we show that this category is fiber-wise algebraically cartesian closed, arithmetical, and strongly protomodular. Moreover, we observe that $\lgrp$ is not action accessible, despite the good behaviour  of centralizers of internal equivalence relations. Finally, we restrict our attention to the subvariety $\lab$ of lattice-ordered abelian groups, showing that it is algebraically coherent; this provides an example of an algebraically coherent category which is not action accessible.

\end{abstract}

\section{Introduction}

A \emph{lattice-ordered group} is a set endowed with both a group structure and a lattice structure such that the underlying order relation is invariant under translations. In other words, a lattice-ordered group can be defined as an algebraic structure of signature $\{ \cdot, e, ^{-1}, \lor, \land \}$ satisfying the axioms of groups, the axioms of lattices, and the axioms related to the distributivity of the group product over both the lattice operations. Therefore, the category of lattice-ordered groups (denoted by $\lgrp$) can be presented as the variety of models associated with the equational theory just described.\\
Recently, lattice-ordered groups have emerged in many areas of mathematics. For instance, in the study of many-valued logic (as shown in \cite{mundici1986interpretation}, the category of lattice-ordered abelian groups with a distinguished order-unit is equivalent to
the one of MV-algebras, which provides algebraic semantics for Łukasiewicz many-valued propositional logic \cite{cignoli2013algebraic}), in the theory of Bézout domains, in complex intuitionistic fuzzy soft set theory, and in varietal questions in universal algebra.\\
Although the notion of lattice-ordered groups is as natural as that of rings or partially ordered groups (it suffices to say that examples of lattice-ordered groups include the set of integers $\mathbb{Z}$, the set of rational numbers $\mathbb{Q}$, and the set of real numbers $\mathbb{R}$ with the usual group sum and the usual order structure), there are currently no studies about this variety from a categorical point of view. The purpose of this work is precisely to explore these aspects.\\
A first observation is that the category of lattice-ordered groups is semi-abelian. In a similar way to how abelian categories describe the properties of the categories of abelian groups and of modules over a ring, the notion of semi-abelian category is aimed at capturing the main algebraic properties of the category of groups. Briefly, a \emph{semi-abelian category} \cite{semiabel} is a pointed finitely cocomplete category which is Barr-exact \cite{barr} and protomodular \cite{bourn1991normalization} (protomodularity, in this context, is equivalent to the Split Short Five Lemma). Examples of semi-abelian categories include, for instance, groups, rings without unit, loops, Lie algebras, Heyting semilattices, etc. However, the notion of semi-abelian category is not as efficient in capturing the properties of groups as the one of abelian category is with respect to abelian groups and modules. Therefore, additional categorical-algebraic conditions have been introduced over the years to get closer to a characterization of the structural properties of the category of groups; among these, one can mention representability of actions \cite{borceux2005representability}, algebraic coherence \cite{artalgebricallycoherent}, and strong protomodularity \cite{strongprot}. This paper is aimed at studying which of
these properties hold in the category of lattice-ordered groups.\\
\indent In Section \ref{Preliminaries} we recall some classical facts about lattice-ordered groups and we focus on the notion of semi-direct product in $\lgrp$.\\
\indent In Section \ref{Centralizers} we study the nature of commutators in $\lgrp$ and we show that every subobject admits a centralizer, which coincides with the classical notion of polar; moreover, we prove that $\lgrp$ is algebraically cartesian closed.\\
\indent In Section \ref{Congruence Distributivity and Arithmetical Category} we give an alternative proof of the known fact that $\lgrp$ is arithmetical using the observation that the only internal group object is the trivial one.\\
\indent In Section \ref{Strong Protomodularity} we show that $\lgrp$ is strongly protomodular; this property implies that, among other things, the commutativity of internal equivalence relations in the Smith-Pedicchio sense \cite{smithpedicchio} is equivalent to the commutativity in the Huq sense \cite{huq} of their associated ideals. Moreover, we observe that in $\lgrp$ every internal equivalence relation admits a centralizer and we provide a description of it.\\
\indent Section \ref{Action Accessibility} is devoted to the study of action accessibility, a property related to the existence of centralizers of internal equivalence relations; here we observe that, despite $\lgrp$ is not action accessible, there is a construction of centralizers which is very close to the one developed in \cite{artactionaccessible} for the action accessible category of rings without unit.\\
\indent Section \ref{Fiber-wise Algebraic Cartesian Closedness} is aimed at proving that the category of lattice-ordered groups is fiber-wise algebraically cartesian closed (i.e.\ each category of points in $\lgrp$ is algebraically cartesian closed); in detail, we show that in the categories of points in $\lgrp$ every subobject admits a centralizer, and we provide a description of it.\\
\indent In Section \ref{Normality of the Higgins Commutator} we study the properties of the Higgins commutator in $\lgrp$; in particular, we prove that $\lgrp$ satisfies the condition of normality of the Higgins commutators showing that the Huq commutator of a pair of ideals (i.e.\ kernels of some arrows) is nothing more than the intersection of the two subobjects.\\
\indent Finally, in Section \ref{Algebraic Coherence}, we focus our attention on the study of the categorical-algebraic properties of the variety of lattice-ordered abelian groups (denoted by $\lab$). To be precise, we show that $\lab$ is algebraically coherent; this condition implies several algebraic properties (such as, for example, strong protomodularity, normality of the Higgins commutator, and the so-called ``Smith is Huq'' condition). Furthermore, we observe that the category of lattice-ordered abelian groups provides an example of an algebraically coherent category that is not action accessible (thus partially solving the Open Problem 6.28 presented in \cite{artalgebricallycoherent}).

\section{Preliminaries}\label{Preliminaries}

In this section, we recall the notion of \emph{lattice-ordered group}. Roughly speaking, a lattice-ordered group is a set endowed with a group structure and a lattice structure such that the group operation is distributive with respect to the lattice operations.

\begin{definition} \label{definizione1}

A \emph{lattice-ordered group} is an algebra $(X, \cdot, e,(-)^{-1}, \lor)$ where:
\begin{itemize}
    \item[LG1)] $(X, \cdot, e, (-)^{-1})$ is a group;
    \item[LG2)] $(X, \lor)$ is a semilattice (i.e.\ $\lor$ is a binary, associative, commutative and idempotent operation on $X$);
    \item[LG3)] for every $x,y,z \in X$ the following equalities hold $$\displaylines{
    x \cdot (y \lor z)=(x \cdot y) \lor (x \cdot z) \text{ and}  \cr
    (x \lor y) \cdot z=(x \cdot z) \lor (y \cdot z).\cr
}$$

\end{itemize}

A \emph{morphism} between two lattice-ordered groups $(X, \cdot, e,(-)^{-1}, \lor)$ and \linebreak $(Y, \cdot, e,(-)^{-1}, \lor)$ is a map $f \colon  X \rightarrow Y$ such that $f$ is both a group homomorphism between $(X, \cdot, e,(-)^{-1})$ and $(Y, \cdot, e,(-)^{-1})$ and a semilattice homomorphism between $(X, \lor)$ and $(Y, \lor)$.\\
The category $\lgrp$ is the category whose objects are the lattice-ordered groups and whose arrows are the morphisms between them.

\end{definition}

Many fundamental algebraic structures naturally admit the structure of lattice-ordered groups. In particular, the set of integers $\mathbb{Z}$, the set of rational numbers $\mathbb{Q}$, and the set of real numbers $\mathbb{R}$ with the usual group sum and the usual order structure are lattice-ordered groups. Moreover, given a totally ordered set $\Gamma$ we can provide a lattice-ordered group structure on the set of order automorphisms $\Aut(\Gamma)$: for every $f,g \in \Aut(\Gamma)$ the group product is defined as the composition $f \circ g$, and $(f \lor g)(x) \coloneqq  \max(f(x),g(x))$ for all $x \in \Gamma$. For these and further examples see e.g.\ \cite{librogruppiret} and \cite{librogrupret2}.\\

Given a lattice-ordered group $X$, we typically denote an algebra by its underlying set. In our notation, we will use $xy$ instead of $x \cdot y$ as usual. Furthermore, we will assume that the product operation precedes the lattice operation. Therefore, we will write $x y \lor z$ to mean $(x y) \lor z$.

In the literature, lattice-ordered groups are usually presented as algebras on the set of operations $\{\cdot, e,(-)^{-1}, \lor, \land \}$ satisfying the group axioms, the lattice axioms and the axioms related to the left and right distributivity of the group operation over both lattice operations. However, in this paper, we have preferred a presentation that does not directly involve the meet operation in order to reduce the number of operations to deal with. In fact, starting from Definition \ref{definizione1} it is always possible to define, in a unique way, the meet operation. To show this, we are about to cite and prove several well-known results, that can be found, for example, in \cite{birkhoff1942lattice}, \cite{librogruppiret}, and \cite{librogrupret2}.
Given a lattice-ordered group $(X, \cdot, e, (-)^{-1}, \lor)$, since $(-)^{-1}$ is an isomorphism of groups between $(X, \cdot, e, (-)^{-1}) \text{ and } (X, \cdot^{\text{op}}, e, (-)^{-1})$, it immediately follows that defining $$x \land y \coloneqq (x^{-1} \lor y^{-1})^{-1}$$ makes $(X, \cdot^{\text{op}}, e, (-)^{-1}, \land)$ a lattice-ordered group, which implies that also\linebreak $(X, \cdot, e, (-)^{-1}, \land)$ is a lattice-ordered group. Moreover, $(X, \lor, \land)$ is a distributive lattice. To prove this, we start by observing that the identity
        \begin{equation}\label{eq1}
           a(x \land y)^{-1}b = (ax^{-1}b) \lor (ay^{-1}b) 
        \end{equation}
        holds. This is immediate since $a(x \land y)^{-1}b = a(x^{-1} \lor y^{-1})b = (ax^{-1}b) \lor (ay^{-1}b)$. The identity \eqref{eq1} implies (setting $a = x$ and $b = y$) that
        \begin{equation}\label{eq2}
            x(x \land y)^{-1}y = x \lor y,
        \end{equation}
        \begin{equation}\label{eq3}
            x = (x \lor y)y^{-1}(x \land y).
        \end{equation}
        To prove that $(X, \lor, \land)$ is a lattice, we observe that $a \lor b = b \text{ if and only if } a \land b = a$. In fact, note that if $a \lor b = b$, then $a = (a \lor b)b^{-1}(a \land b) = b(b^{-1}(a \land b)) = a \land b$, and if $a \land b = a$, then $b = (b \lor a)a^{-1}(b \land a) = (b \lor a)a^{-1}a = b \lor a$. Finally, to show the distributivity, we recall the well-known fact that a lattice is distributive if and only if the implication $$[x \land b = y \land b \text{ and } x \lor b = y \lor b] \implies x = y$$ holds. This is true because if $x \land b = y \land b$ and $x \lor b = y \lor b$, then, by the identity \eqref{eq3}, we get $x = (x \lor b)b^{-1}(x \land b) = (y \lor b)b^{-1}(y \land b) = y$.\\
        Moreover, when $\cdot$ is commutative, identity \eqref{eq2} implies $x  y = (x \land y)(x \lor y)$. This fact will be crucial in the last part of this paper.
        
\begin{definition}[\cite{librogruppiret,librogrupret2}]
Let $X$ be an object of $\lgrp$. For every $x \in X$ we define $$x^+ \coloneqq  x \lor e \textit{, } x^- \coloneqq  x \land e \textit{ and } |x| \coloneqq  x \lor \inv{x};$$
$x^+$ is called the \emph{positive part} of $x$, $x^-$ the \emph{negative part} of $x$, and $|x|$ the \emph{absolute value} of $x$.

\end{definition}

The previous definition is useful in order to see that, in a lattice-ordered group, every element can be written as the product of a positive element and a negative one. Specifically, given an object $X$ of $\lgrp$, we define the positive cone of $X$ as $$P \coloneqq  \{ x \in X \, | \, x \geq e \}.$$ It is a well-known fact that $X$ is generated by its positive cone (i.e.\ $X=P \inv{P}$). We offer a proof of this fact by essentially following the approach outlined in \cite{birkhoff1942lattice}, \cite{librogruppiret}, and \cite{librogrupret2}. To do this, we show that 
\begin{equation}\label{eq4}
    x=x^+x^-,
\end{equation}
for every $x \in X$. Clearly, \eqref{eq4} follows by \eqref{eq3} setting $y=e$. Another well-known and interesting fact regarding the behavior of elements in a lattice-ordered group is that
\begin{equation}
    |x| = x^+(x^-)^{-1}.
\end{equation}
We prove this identity, too (once again, the proof follows what is presented in \cite{birkhoff1942lattice}, \cite{librogruppiret}, and \cite{librogrupret2}). Note that $x \lor x^{-1} \geq x \land x^{-1}$, hence $(x \lor x^{-1})(x \land x^{-1})^{-1} \geq e$, and so $(x \lor x^{-1})^2 \geq e$. Now, since $(a \land e)^n = (a^n \land e) \land (a \land e)^{n-1}$ (to see this, it suffices to expand and observe that $(a \land e)^n = a^n \land a^{n-1} \land \dots \land e$), it follows that $a^n \geq e$ implies $(a \land e)^n=(a \land e)^{n-1}$, and so $a \geq e$. Therefore, $x \lor x^{-1} \geq e$. To conclude, we compute $x^+(x^-)^{-1} = (x \lor e)(x^{-1} \lor e) = e \lor x \lor x^{-1} \lor e = x \lor x^{-1}$.\\
Furthermore, the notion of positive part is extremely useful in order to characterize group homomorphisms between lattice-ordered groups which are, in addition, morphisms of $\lgrp$. In fact, the following holds:

\begin{lemma} \label{truccomorfismi}

Let $X,Y$ be two objects of $\lgrp$. A map $f \colon X \rightarrow Y$ is a morphism of $\lgrp$ if and only if $f$ preserves the group product and $$f(x \lor e)=f(x) \lor e \text{, for all x} \in X.$$

\begin{proof}

One implication is trivial. So, let us suppose that $f$ preserves the group product (i.e.\ is a group homomorphism) and $f(x \lor e)=f(x) \lor e \text{ for every }x \in X.$ We have to prove that $$f(x \lor y)=f(x) \lor f(y) \text{, for all } x,y \in X.$$ We have $x \lor y=(x \inv{y} \lor e)y$, hence $f(x \lor y)=f(x \inv{y} \lor e)f(y)$ and, by assumption,
$$f(x \inv{y} \lor e)f(y)=(f(x \inv{y}) \lor e)f(y)=(f(x)\inv{f(y)}\lor e)f(y)=f(x) \lor f(y). \eqno{\qedhere}$$
\end{proof}

\end{lemma}

This last result follows essentially from Theorem 9 of \cite{birkhoff1942lattice}.\\

Now, we want to provide a description of the ideals in the variety $\lgrp$. In a category where it makes sense to speak of a kernel of a morphism (for example a pointed finitely complete category), a subobject of $X$ is called an \emph{ideal} (or \emph{normal subobject}) if it is the kernel of some morphism. A detailed study of the notion of an ideal in the variety $\lgrp$ can be found in \cite{librogruppiret} and \cite{librogrupret2}.\\
First of all, we have to recall the definition of a \emph{convex} subset. Given an object $X$ of $\lgrp$ and a subset $S \subseteq X$, $S$ is said to be convex if for every $a,b \in S$ and every $x \in X$, if $a \leq x \leq b$ then $x \in S$.\\ A subobject $A \leq X$ is an ideal if and only if it is normal (in the classical sense) as a subgroup and it is a convex subset.\\

The aim of the following proposition is to describe the notion of convexity only with terms. This characterization will be crucial for the purpose of working with semi-direct products.

\begin{proposition}

Let $X$ be an object of $\lgrp$ and $A \leq X$ a subalgebra. $A$ is convex if and only if for every $a_1,a_2 \in A$ and $x, y \in X$ one has $$(a_1x \lor a_2y)\inv{(x \lor y)} \in A.$$

\begin{proof}

Let us suppose that $A$ is convex. We consider the following inequalities: $$((a_1 \land a_2)x) \lor ((a_1 \land a_2)y) \leq a_1x \lor a_2y \leq ((a_1 \lor a_2)x)\lor((a_1 \lor a_2)y),$$
hence $$(a_1 \land a_2)(x \lor y) \leq a_1x \lor a_2y \leq (a_1 \lor a_2)(x \lor y),$$
and so
$$(a_1 \land a_2) \leq (a_1x \lor a_2y)\inv{(x \lor y)} \leq (a_1 \lor a_2).$$
Thus, since $A$ is convex, we deduce $(a_1x \lor a_2y)\inv{(x \lor y)} \in A$.\\
Conversely, let us suppose that for every $a_1,a_2 \in A$ and $x, y \in X$ one has $(a_1x \lor a_2y)\inv{(x \lor y)} \in A$. Let us take an element $x \in X$ and two elements $a_1, a_2 \in A$ such that $a_1 \leq x \leq a_2$. We want to prove that $x$ belongs to $A$, We observe that $a_1 \lor x=x$ and $a_2 \land x= x$. Therefore,
\begin{align*}
    x&=(a_1 \lor x)\inv{x}(x \land a_2)=(a_1 \lor x)(e \land \inv{x}a_2)=(a_1 \lor x)\inv{(e \lor \inv{a_2}x)}\\
    &=(a_1 e \lor a_2 \inv{a_2}x)\inv{(e \lor \inv{a_2}x)} \in A;
\end{align*}
where the last term belongs to $A$ by assumption.
\end{proof}

\end{proposition}

This means that a subobject $A \leq X$ in $\lgrp$ is normal if and only if for every $a_1,a_2,a \in A$ and $x, y,z \in X$ one has $(a_1x \lor a_2y)\inv{(x \lor y)} \in A$ and $\conj{z}{a} \in A$.\\

In the following part of this section, we will deal with the notion of a semi-abelian category. The concept of a semi-abelian category aims to capture some of the common algebraic properties of the category of groups; among the examples of semi-abelian category we can find those of groups, rings without unit, Lie algebras, and Heyting semilattices. 

\begin{definition}[\cite{semiabel}] A pointed category (i.e.\ a category with a zero object) $\mathbb{C}$ is \emph{semi-abelian} if:
\begin{itemize}
    \item it is Barr-exact \emph{\cite{barr}} (which means that $\mathbb{C}$ is a regular category in which every internal equivalence relation is a kernel pair);
    \item it has finite coproducts;
    \item it is protomodular \emph{\cite{bourn1991normalization}} (in this context, this is equivalent to the Split Short Five Lemma holding in $\mathbb{C}$).
\end{itemize}

\end{definition}

In \cite{varprot}, Theorem 1.1 the authors provided, in the case of a variety $\mathbb{V}$ of universal algebras, a characterization for protomodularity depending on terms. In fact, the authors proved that a variety $\mathbb{V}$ is protomodular if and only if it has $0$-ary terms $e_1, \dots, e_n$, binary terms $t_1, \dots, t_n$ and an $(n+1)$-ary term $t$ satisfying the identities $$t(x,t_1(x,y), \dots ,t_n(x,y))=y \text{ and } t_i(x,x)=e_i$$ for all $i=1, \dots ,n$.\\

Since every variety of $\Omega$-groups is semi-abelian (see Example 2.6 of \cite{semiabel}) we obtain:

\begin{proposition}

$\lgrp$ is a semi-abelian category.

\end{proposition}

As shown in \cite{bourn1998protomodularity}, in every semi-abelian category there exist semi-direct products in a categorical sense. In the category of groups, the categorical semi-direct product coincides with the classical one. Now we can describe semi-direct products in the category $\lgrp$. In order to do this, we will apply, in the next proposition, the results provided in \cite{artprodottisemidirettialgebre}.

\begin{proposition} \label{semidiretto}

Let $p \colon A \rightarrow B$ be a split epimorphism in $\lgrp$ with fixed section $s \colon  B \rightarrow A$, and $k \colon K \rightarrow A$ a kernel of $p$. Without loss of generality let us suppose that $K,B$ are subalgebras of $A$ and $k,s$ are the inclusions of subalgebras. Then $A$ is isomorphic (as a lattice-ordered group) to the set $K \times B$ endowed with the operations
\begin{itemize}
    \item $(k_1,b_1)(k_2,b_2)=(k_1 b_1 k_2 b_1^{-1},b_1b_2)$,
    \item $(k_1,b_1)\lor(k_2,b_2)=((k_1b_1 \lor k_2b_2)\inv{(b_1 \lor b_2)},b_1 \lor b_2)$,
\end{itemize}
(which takes the name of \emph{semi-direct product} and will be denoted by $K \rtimes B$)
via the morphism 
\begin{align*}
\varphi \colon  & K \rtimes B \rightarrow A \\
        & (k,b) \mapsto kb.
\end{align*}
Moreover, considering the following diagram in $\lgrp$:
\[\begin{tikzcd}
	K & {K \rtimes B} & B \\
	K & A & B
	\arrow["k", from=2-1, to=2-2]
	\arrow["p", shift left=1, from=2-2, to=2-3]
	\arrow["s", shift left=1, from=2-3, to=2-2]
	\arrow[equal, from=1-1, to=2-1]
	\arrow[equal, from=1-3, to=2-3]
	\arrow["{p_B}", shift left=1, from=1-2, to=1-3]
	\arrow["{i_B}", shift left=1, from=1-3, to=1-2]
	\arrow["{i_K}", from=1-1, to=1-2]
	\arrow["\varphi"', from=1-2, to=2-2]
\end{tikzcd}\]
where $i_K(k)=(k,e)$, $i_B(b)=(e,b)$ and $p_B(k,b)=b$, we have $\varphi i_K=k,p \varphi=p_B$, and $\varphi i_B=s$.
\end{proposition}

\section{Centralizers and Algebraic Cartesian Closedness}\label{Centralizers}

In this section we study, from a categorical point of view, the commutativity of subobjects in the variety $\lgrp$.\\ In order to introduce the topic, we mention some known results related to the category of groups. Given a group $G$ and two subgroups $A,B \leq G$, the condition that, for every $a \in A$ and $b \in B$, $ab=ba$ can be reformulated in the following equivalent way: there exists a group homomorphism $\varphi \colon  A \times B \rightarrow G$ making the following diagram commutative:
\[\begin{tikzcd}
	A && B \\
	& {A \times B} \\
	\\
	& G.
	\arrow["{(id_A,0)}"{description}, from=1-1, to=2-2]
	\arrow["{(0,id_B)}"{description}, from=1-3, to=2-2]
	\arrow[bend right, hook, from=1-1, to=4-2]
	\arrow[bend left, hook', from=1-3, to=4-2]
	\arrow["\varphi", from=2-2, to=4-2]
\end{tikzcd}\]
Moreover, it is easy to show that $\varphi$ must be the group product and, therefore, it is necessarily unique. Hence, with the aim of generalising the notion of commutativity, we must place ourselves in a context in which a morphism  $\varphi$ of this type is unique. This reasoning justifies the following definition:

\begin{definition}[\cite{unitalcat}] A pointed category $\mathbb{C}$ with finite products is \emph{unital} if, for $X$ and $Y$ objects of $\cat{C}$, the pair of morphisms $(id_X,0) \colon X \rightarrow X \times Y$, $(0,id_Y) \colon Y \rightarrow X \times Y$ is jointly extremally epimorphic.

\end{definition}

To be more explicit, a pair of arrows $f \colon  A \rightarrow B$ and $g \colon  C \rightarrow B$ of a category $\mathbb{C}$ is said to be jointly extremally epimorphic when, for every commutative diagram

\[\begin{tikzcd}
	& M \\
	A & B & C,
	\arrow["f"', from=2-1, to=2-2]
	\arrow["g", from=2-3, to=2-2]
	\arrow["{f'}", from=2-1, to=1-2]
	\arrow["{g'}"', from=2-3, to=1-2]
	\arrow["m", from=1-2, to=2-2]
\end{tikzcd}\]

if $m$ is a monomorphism, then $m$ is an isomorphism.\\

It has been shown in \cite{semi} that every semi-abelian category is unital.\\

We are ready to mention the generalized notion of commutativity between subobjects.

\begin{definition}[\cite{huq}]
Let $\cat{C}$ be a unital category. Two subobjects $a \colon  A \mono X$ and $b \colon  B \mono X$ of $X$ are said to \emph{cooperate} (or \emph{commute in the sense of Huq}, and we write $[a,b]=0$) if there exists a (necessarily unique) morphism $\varphi \colon  A \times B \rightarrow X$ (called \emph{cooperator}) such that the following diagram commutes:
\[\begin{tikzcd}
	A && B \\
	& {A \times B} \\
	\\
	& X.
	\arrow[ bend right, tail , "{a}"', from=1-1, to=4-2]
	\arrow[ bend left, tail , "{b}", from=1-3, to=4-2]
	\arrow["{(id_A,0)}"{description}, from=1-1, to=2-2]
	\arrow["{(0,id_B)}"{description}, from=1-3, to=2-2]
	\arrow["{\exists  \varphi}", from=2-2, to=4-2]
\end{tikzcd}\]
Given a subobject $a \colon  A \mono X$, the \emph{centralizer} of $a$ in $X$, if it exists, is the greatest subobject of $X$ that cooperates with $a$.
\end{definition}

Now, let us recall the definition of \emph{orthogonal} subobjects of a lattice-ordered group. This concept will be essential in order to study the condition of cooperation.

\begin{definition}[\cite{birkhoff1942lattice}] Let $X$ be an object of $\lgrp$. Two elements $a,b \in X$ are called \emph{orthogonal} if $$|a| \land |b|=e.$$
Two subsets $A,B \subseteq X$ are called \emph{orthogonal} (and one writes $A \perp B$) if, for every $a \in A$ and for every $b \in B$, $a$ and $b$ are orthogonal as elements.

\end{definition}

It is a known fact that two orthogonal subobjects of a lattice-ordered group commute as subgroups. More generally, if $a$ and $b$ are orthogonal then $ab=ba$. A proof of this can be found, for instance, in Proposition 2.2.10 of \cite{librogruppiret}.

\begin{proposition} \label{central}

Let $X$ be an object of $\lgrp$ and $A,B \leq X$ two subobjects. Then $A$ and $B$ cooperate if and only if $A \perp B$.

\begin{proof}

($\Rightarrow$) The cooperator $\varphi \colon  A \times B \rightarrow X$ is given by $\varphi(a,b)=a  b$. In fact, since $\varphi$ preserves the group operation, we have $$\varphi(a,b)=\varphi(a,e)\varphi(e,b)=\varphi(i_A(a))\varphi(i_B(b))=ab.$$ 
We observe that, for every $a \in A$ and $b \in B$, $(|a|,e) \land (e,|b|)=(e,e)$ holds. So, since $\varphi$ preserves the lattice operations, we get $$e=\varphi(e,e)=\varphi((|a|,e) \land (e,|b|))=\varphi(|a|,e) \land \varphi(e,|b|)=|a| \land |b|.$$
($\Leftarrow$) If a cooperator $\varphi$ exists then it must be the group multiplication because of what we observed at the beginning of the proof. In fact, in order to guarantee the existence of a cooperator, we only have to prove 
\begin{equation}\label{coppia di eq}
    ab=ba \text{ and } (ab)^+=a^+b^+
\end{equation}
for all $a \in A$ and $b \in B$. Indeed, if $ab=ba$ for every $a \in A$ and $b \in B$, then $\varphi((a,b)(c,d))=\varphi(ac,bd)=acbd=abcd=\varphi(a,b) \varphi(c,d)$ since, by assumption, $cb=bc$; furthermore, if $(ab)^+=a^+b^+$ for every $a \in A$ and $b \in B$, then $\varphi((a,b)^+)=\varphi((a,b) \lor (e,e))=\varphi(a^+,b^+)=a^+b^+=(ab)^+=\varphi(a,b)^+=\varphi(a,b) \lor e$, and thus we can apply Lemma \ref{truccomorfismi} to say that $\varphi$ is a morphism of lattice-ordered groups. As recalled above, the first equality of \eqref{coppia di eq} holds since $A \perp B$. For the second one, we observe that $a^+b^+=(a \lor e)(b \lor e)=ab \lor a \lor b \lor e$ and $(ab)^+=ab \lor e$; so we have to prove $a \lor b \leq ab \lor e$. Since $|a| \land |b|=e$, we have $(a \lor \inv{a}) \land (b \lor \inv{b})=e$ and, by distributivity, we get $$(a \land b)\lor(a \land \inv{b})\lor(\inv{a} \land b)\lor(\inv{a} \land \inv{b})=e.$$
Hence, $\inv{a} \land b \leq e$ implies $a \lor \inv{b} \geq e$ and, multiplying by $b$ on the right, we obtain $ab \lor e \geq b$; with a similar argument we get $ab \lor e \geq a$. Finally, considering the last two inequalities, we conclude that $ab \lor e \geq a \lor b$.
\end{proof}

\end{proposition}

As a direct consequence of this proof we have the following:

\begin{corollary}\label{magie tra ortogonali}

Let $X$ be an object of $\lgrp$ and $A,B \leq X$ two subobjects. Then $A \perp B$ if and only if, for every $a \in A$ and $b \in B$, the following equalities hold: $$ab=ba \textit{ and } (ab)^+=a^+b^+.$$
    
\end{corollary}
The previous result will be of fundamental importance in the next sections.\\

We recall the notion of \emph{polar} of a subset $S$ of a lattice-ordered group (i.e.\ the set of elements orthogonal to each element of $S$). We will show that the polar of a subobject is nothing but the centralizer of the subobject. Hence, we will exhibit some properties of centralizers related to being ideals.

\begin{proposition}[\cite{librogruppiret}, Proposition 1.2.6] Let $X$ be an object of $\lgrp$ and $S \subseteq X$ a non-empty subset. Then the set $\pol{S} \coloneqq  \{x \in X \mid \textit{for each } s \in S \textit{ } |x| \land |s|=e \}$ (called the \emph{polar} of $S$) is a convex subalgebra of $X$.

\end{proposition}

\begin{lemma} \label{coniugiopolari}

Let $X$ be an object of $\lgrp$ and $S \subseteq X$ a non-empty subset of $X$ closed under conjugation. Then $\pol{S}$ is an ideal of $X$.

\begin{proof}

First of all we observe that, for all $x,y \in X$, one has $$|\conj{x}{y}|=\conj{x}{y} \lor \conj{x}{\inv{y}}=\inv{x}(y \lor \inv{y})x=\conj{x}{|y|}.$$
We want to show that $|\conj{x}{y}| \land |s|=e$, for every $x \in X$, $y \in \pol{S}$, and $s \in S$. We observe that $|\conj{x}{y}| \land |s|=\conj{x}{|y|} \land |s|=\conj{x}{(|y| \land \jnoc{x}{|s|})}=\conj{x}{(|y| \land |\jnoc{x}{s}|)}=\conj{x}{(|y| \land |\overline{s}|)}=e$, where $\jnoc{x}{s}=\overline{s} \in S$ because $S$ is closed under conjugation, and $|y| \land |\overline{s}|=e$ (since $y \in \pol{S}$).
\end{proof}

\end{lemma}

\begin{corollary}

Let $X$ be an object of $\lgrp$ and $A,B \leq X$ two subobjects. $A$ and $B$ cooperate if and only if $B \subseteq \pol{A}$. Therefore, $\pol{A} \leq X$ is the centralizer of $A \leq X$. Moreover, if $A$ is an ideal of $X$, then $\pol{A}$ is an ideal of $X$, too.

\end{corollary}

Finally, we recall a property that is strictly related to the existence of centralizers. It is well known that a category $\mathbb{E}$ with finite products is cartesian closed if and only if for every object $Y$ of $\mathbb{E}$ the change-of-base functor $\tau_Y^* \colon  \mathbb{E} \rightarrow \mathbb{E}/Y$ along the terminal arrow $\tau_Y \colon  Y \rightarrow \terminal$ has a right adjoint. For algebraic categories, such adjoints rarely exist, but it turns out to be of interest to consider a variation of this notion; this leads to:

\begin{definition}[\cite{acc}]

A category $\mathbb{C}$ is \emph{algebraically cartesian closed (a.c.c.)} if for every object $X$ of $\mathbb{C}$ the change-of-base functor $\tau_X^* \colon  \Pt{\terminal}{\mathbb{C}} \rightarrow \Pt{X}{\mathbb{C}}$ has a right adjoint, where $\tau_X \colon X \rightarrow \terminal$ is the unique arrow from $X$ to the terminal object.

\end{definition}

In \cite{acc} the authors show that the existence of such adjoints is related to the existence
of cofree structures for the split epimorphisms $p_Y \colon  Y \times X \rightarrow Y$ in $\pt{Y}{\mathbb{E}}$ with fixed section $(id_Y,u)$, where $u \colon  Y \rightarrow X$ can be chosen to be a monomorphism. 

\begin{proposition}[\cite{acc}, Proposition 1.2] A unital category $\mathbb{C}$ is algebraically cartesian closed if and only if, for every $X$ object of $\mathbb{C}$, each subobject of $X$ has a centralizer.

\end{proposition}

\begin{corollary}

The category $\lgrp$ is algebraically cartesian closed.

\end{corollary}

\section{Congruence Distributivity and Arithmeticity} \label{Congruence Distributivity and Arithmetical Category}

It is a widely known fact that in the category $\lgrp$ the lattice of congruences on any object is distributive (see Theorem 21 of \cite{birkhoff1942lattice}). In this section we provide an alternative proof of this fact based on categorical tools. We recall that a category is a \emph{Mal'tsev category} \cite{carboni1991diagram} if it is finitely complete and if every internal reflexive relation is an internal equivalence relation. If the category is regular, this notion is equivalent to the following property: for every object $X$ and for every pair of internal equivalence relations $(s_1,s_2) \colon  S \mono X \times X$ and $(r_1,r_2) \colon  R \mono X \times X$ one has $R \circ S=S \circ R$; in detail,  $R \circ S \colon  \mono X \times X$ is defined as the regular image of $(p_1,p_3)$, where $(p_1,p_3)$ is given by the following diagram:
\[\begin{tikzcd}
	{R \times_XS} & S & X \\
	R & X \\
	X.
	\arrow["{\pi_S}"', from=1-1, to=1-2]
	\arrow["{\pi_R}", from=1-1, to=2-1]
	\arrow["{r_2}"', from=2-1, to=2-2]
	\arrow["{s_1}", from=1-2, to=2-2]
	\arrow["{r_1}", from=2-1, to=3-1]
	\arrow["{s_2}"', from=1-2, to=1-3]
	\arrow["\lrcorner"{anchor=center, pos=0.125}, draw=none, from=1-1, to=2-2]
	\arrow[bend right, "{p_1}"{description}, from=1-1, to=3-1]
	\arrow[bend left, "{p_3}"{description}, from=1-1, to=1-3]
\end{tikzcd}\]
The composite $S \circ R$ can be defined in a similar way. Moreover, if the category is a variety of universal algebras, the property of being a Mal'tsev category is equivalent to the existence of a ternary term $p(x,y,z)$ (called \emph{Mal'tsev term}) satisfying the axioms $$p(x,x,z)=z \text{ and } p(x,y,y)=x,$$ for every object $X$ and for every $x,y,z \in X$. Therefore, if the theory contains a group operation, the associated variety is a Mal'tsev category: in fact, a Mal'tsev ternary term is $p(x,y,z) \coloneqq  x \inv{y} z$.\\
Then, we immediately get the following result:

\begin{corollary}

The category $\lgrp$ is a Mal'tsev category.

\end{corollary}

If $\mathbb{C}$ is a Barr-exact category with coequalizers then, for every object $X$ of $\mathbb{C}$, the set $\Eq(X)$ of internal equivalence relations on $X$ is a lattice; given two internal equivalence relations $(s_1,s_2) \colon  S \mono X \times X$ and $(r_1,r_2) \colon  R \mono X \times X$, the meet $S \land R$ is defined as the meet of subobjects of $X \times X$, and the join $(t_1,t_2) \colon  S \lor R \mono X \times X$ is defined as the kernel pair of $q=\coeq(v_1,v_2)$, where $(v_1,v_2) \colon  V \mono X \times X$ is the join of $S$ and $R$ as subobjects of $X \times X$ (we recall that the join, as subobjects, of two internal equivalence relations is not, in general, an internal equivalence relation). Thanks to the previous observations, the classical notion of arithmetical variety of universal algebras can be extended to a categorical context as follows:

\begin{definition}[\cite{aritmetica}] A Barr-exact category with coequalizers $\mathbb{C}$ is \emph{arithmetical} if it is a Mal'tsev category and, for any object $X$ of $\mathbb{C}$, the lattice $\Eq(X)$ of internal equivalence relations on $X$ is distributive. 

\end{definition}

It is a known fact (a proof of this can be found in \cite{semi}) that the property of being an arithmetical category is related to the absence of non-trivial \emph{internal group objects} in the category. In fact, the following holds:

\begin{proposition}[\cite{semi}, Proposition 2.9.9] Let $\mathbb{C}$ be a semi-abelian category. If in $\mathbb{C}$ the only internal group object is the zero object then $\mathbb{C}$ is arithmetical.

\end{proposition}

\begin{proposition}

The only internal group object in $\lgrp$ is the zero object.

\begin{proof}

Given an internal group $X$ in $\lgrp$, with multiplication $\mu  \colon  X \times X \rightarrow X$ and neutral element $\eta  \colon  \{* \} \rightarrow X$, we want to show that $X = \{e\}$.
It is not difficult to see that $\eta(*)=e$ and $\mu (x,y)=xy$ for all $x,y \in X$. This determines a cooperator between $X$ and itself. Hence, we get $X \perp X$ and so $x^+=x^+ \land x^+=|x^+| \land |x^+|=e$, for every $x \in X$. With a similar argument we can prove $x^-=e$, for every $x \in X$. Finally, recalling that $x=x^+x^-$, we obtain $x=e$ for every $x \in X$ (i.e.\ $X=\{ e \}$).

\end{proof}

\end{proposition}

\begin{corollary}

$\lgrp$ is an arithmetical category.

\end{corollary}

\section{Strong Protomodularity} \label{Strong Protomodularity}

Given a category $\mathbb{C}$, we denote by $\point (\mathbb{C})$ the category whose objects are the diagrams in $\mathbb{C}$ of the form
\[\begin{tikzcd}
	A & B
	\arrow["p", shift left=1, from=1-1, to=1-2]
	\arrow["s", shift left=1, from=1-2, to=1-1]
\end{tikzcd}\]
where $ps=id_B$, and whose arrows are the pairs $(f,g)$  of arrows of $\mathbb{C}$
\[\begin{tikzcd}
	A & B \\
	C & D
	\arrow["p", shift left=1, from=1-1, to=1-2]
	\arrow["s", shift left=1, from=1-2, to=1-1]
	\arrow["q", shift left=1, from=2-1, to=2-2]
	\arrow["r", shift left=1, from=2-2, to=2-1]
	\arrow["f"', shift right=1, from=1-1, to=2-1]
	\arrow["g", from=1-2, to=2-2]
\end{tikzcd}\]
such that $qf=gp$ and $fs=rg$. We denote by $\pi \colon  \point (\mathbb{C}) \rightarrow \mathbb{C}$ the functor that associates to every split extension (i.e.\ an object of $\point (\mathbb{C})$) $(p,s)$ the codomain of $p$.

\begin{definition}[\cite{strongprot}] A finitely complete category $\mathbb{C}$ is \emph{strongly protomodular} when all the change-of-base functors of $\pi \colon  \point (\mathbb{C}) \rightarrow \mathbb{C}$ reflect both isomorphisms and normal monomorphisms (in the semi-abelian case a monomorphism is normal if and only if it is the kernel of some arrow).

\end{definition}

In \cite{strongprot} the author shows that, if $\mathbb{C}$ is a pointed protomodular category, there is a characterization of strong protomodularity related to the stability of kernels. Let us consider a diagram in $\mathbb{C}$ of the form
\[\begin{tikzcd}
	X & A & B \\
	Y & C & B
	\arrow["l"', from=2-1, to=2-2]
	\arrow["m"', from=1-1, to=2-1]
	\arrow["k", from=1-1, to=1-2]
	\arrow["f"', from=1-2, to=2-2]
	\arrow["p"', shift right=1, from=1-2, to=1-3]
	\arrow["s"', shift right=1, from=1-3, to=1-2]
	\arrow["q"', shift right=1, from=2-2, to=2-3]
	\arrow["r"', shift right=1, from=2-3, to=2-2]
	\arrow[equal, from=1-3, to=2-3]
\end{tikzcd}\]

 where $k=\kerr(p)$, $ps=id_B$, $l=\kerr(q)$, $qr=id_B$, $m$ is a normal monomorphism and the right-rightward square, the right-leftward square, and the left square commute. Then $\mathbb{C}$ is strongly protomodular if and only if the composite $lm$ is a normal monomorphism, for every diagram of this form.

\begin{proposition}

$\lgrp$ is a strongly protomodular category.

\begin{proof}

Let us consider the following commutative diagram (without loss of generality we can assume that the monomorphisms are inclusions):

\[\begin{tikzcd}
	X & A & B \\
	Y & C & B
	\arrow[hook, from=2-1, to=2-2]
	\arrow[hook, from=1-1, to=2-1]
	\arrow[hook, from=1-1, to=1-2]
	\arrow["f"', from=1-2, to=2-2]
	\arrow["p"', shift right=1, from=1-2, to=1-3]
	\arrow["s"', shift right=1, from=1-3, to=1-2]
	\arrow["q"', shift right=1, from=2-2, to=2-3]
	\arrow["r"', shift right=1, from=2-3, to=2-2]
	\arrow[equal, from=1-3, to=2-3]
\end{tikzcd}\]

 where $ps=id_B$, $qr=id_B$, $X$ is the ideal of $A$ determined by $\kerr(p)$, $Y$ is the ideal of $C$ determined by $\kerr(q)$, $X$ is an ideal of $Y$ and the right-rightward square, the right-leftward square, and the left square commute. We want to show that $X$ is an ideal of $C$, too.
\begin{itemize}
    \item $X$ is a normal subgroup of $C$. This immediately follows from the fact that the category of groups is strongly protomodular, as established in \cite{strongprot}.
    \item $X$ is a convex subset of $C$. We know that, for every $y \in Y$, if $x_1 \leq y \leq x_2$, with $x_1,x_2 \in X$, then $y \in X$ and, for every $c \in C$, if $y_1 \leq c \leq y_2$, with $y_1, y_2 \in Y$, then $c \in Y$. So, given an element $c \in C$ such that $x_1 \leq c \leq x_2$, with $x_1,x_2 \in X$, we have $c \in Y$ (since $x_1,x_2 \in Y$) and thus $c \in X$. \qedhere
\end{itemize}
\end{proof}

\end{proposition}

In the final part of this section we study, in the case of $\lgrp$, the consequences of strong protomodularity relatively to the commutativity, in the Smith-Pedicchio sense, of internal equivalence relations. In particular, we show that every internal equivalence relation admits a centralizer. Let us begin by recalling the necessary notions to deal with this subject.

\begin{definition}[\cite{smithpedicchio}, \cite{bourn2002centrality}]

Let $\mathbb{C}$ be a Mal’tsev category and

\[\begin{tikzcd}
	R & {X,} & S & X
	\arrow["{r_1}", shift left=3, from=1-1, to=1-2]
	\arrow["{r_2}"', shift right=3, from=1-1, to=1-2]
	\arrow["{\delta_R}"{description}, from=1-2, to=1-1]
	\arrow["{s_1}", shift left=3, from=1-3, to=1-4]
	\arrow["{s_2}"', shift right=3, from=1-3, to=1-4]
	\arrow["{\delta_S}"{description}, from=1-4, to=1-3]
\end{tikzcd}\]

 a pair of internal equivalence relations on an object $X$ of $\mathbb{C}$. We say that $(R,r_1,r_2)$ and $(S,s_1,s_2)$ \emph{commute in the Smith-Pedicchio sense} (and we write $[R,S]=0$) if, given the following diagram:

\[\begin{tikzcd}
	{R \times_XS} & S \\
	R & X
	\arrow["{s_1}"', shift right=1, from=1-2, to=2-2]
	\arrow["{r_2}", shift left=1, from=2-1, to=2-2]
	\arrow["{\delta_S}"', shift right=1, from=2-2, to=1-2]
	\arrow["{\delta_R}", shift left=1, from=2-2, to=2-1]
	\arrow["{\pi_R}"', shift right=1, from=1-1, to=2-1]
	\arrow["{\pi_S}", shift left=1, from=1-1, to=1-2]
	\arrow["{\tau_R}"', shift right=1, from=2-1, to=1-1]
	\arrow["{\tau_S}", shift left=1, from=1-2, to=1-1]
	\arrow["\lrcorner"{anchor=center, pos=0.125}, draw=none, from=1-1, to=2-2]
\end{tikzcd}\]

 where $R \times_X S$ is the pullback of $r_2$ through $s_1$, $\tau_R=(id_R, \delta_S r_2)$ and $\tau_S=(\delta_R s_1,id_S)$ are induced by the universal property, there exists a unique morphism $p \colon R \times_X S \rightarrow X$ (called \emph{connector} between $R$ and $S$) such that $p \tau_S=s_2$ and $p \tau_R=r_1$. The \emph{centralizer} of an internal equivalence relation $(R,r_1,r_2)$ on $X$, if it exists, is the
greatest internal equivalence relation on $X$ which commutes with $(R,r_1,r_2)$.

\end{definition}

It has been shown in Proposition 3.2 of \cite{strongprotSH} that, in a pointed Mal’tsev
category, if two internal equivalence relations $(R,r_1,r_2)$
and $(S,s_1,s_2)$ commute in the Smith-Pedicchio sense, then necessarily their associated normal subobjects $j_R$ and $j_S$ commute in the Huq sense, where $j_R \coloneqq \kerr(q_R)$ and $j_S \coloneqq \kerr(q_S)$, with $q_R \coloneqq \coeq(r_1,r_2)$ and $q_S \coloneqq \coeq(s_1,s_2)$. Briefly, $[R,S]=0$ implies $[j_R,j_S]=0$. The converse is not true, in general. We say that a pointed Mal'tsev category satisfies the so-called \emph{Smith is Huq} condition (SH) if $[j_R,j_S]=0$ implies $[R,S]=0$. It has been proved in Theorem 6.1 of \cite{strongprotSH} that in every pointed strongly protomodular category the Smith is Huq condition holds.\\

In \cite{smith2006mal}, page 39, it is shown that, in every Mal'tsev variety, internal equivalence relations have centralizers. We provide an explicit description of these centralizers in $\lgrp$.\\
Since $\lgrp$ is a semi-abelian category, we have, for every object $X$, an order-preserving bijection $\varphi$ between $\Eq(X)$ and the lattice $\Ideals(X)$ of ideals of $X$, where $\varphi(R) \coloneqq I_R$, with $I_R \coloneqq \{ x \in X \, | \, (x,e) \in R \}$. Given two internal equivalence relations $R \leq X \times X$ and $S \leq X \times X$, we know that $[R,S]=0$ if and only if $[\varphi(R), \varphi(S)]=0$ ($\lgrp$ is a strongly protomodular category, hence (SH) holds). Moreover, given an internal equivalence relation $R$ on $X$, we recall from Lemma \ref{coniugiopolari} that the centralizer $\pol{I_R} \leq X$ of the ideal $I_R$ associated with $R$ is an ideal. Therefore, since $\varphi$ is an order-preserving bijection and $\lgrp$ satisfies (SH), the centralizer of $R$ in $X$ is $\inv{\varphi}(\pol{I_R})$.

\section{Action Accessibility} \label{Action Accessibility}

To approach the topic covered in this section, let us recall a fundamental concept in the category $\grp$ of groups: the notion of \emph{split extension}.

\begin{definition}

Let $\mathbb{C}$ be a pointed protomodular category. A \emph{split extension} of $\mathbb{C}$ is a diagram of the form
\[\begin{tikzcd}
	X & A & B,
	\arrow["k", from=1-1, to=1-2]
	\arrow["s"', shift right=1, from=1-3, to=1-2]
	\arrow["p"', shift right=1, from=1-2, to=1-3]
\end{tikzcd}\]
where $k=\ker(p)$ and $ps=id_B$.\\
We denote by $\splt{X}{\mathbb{C}}$ the category whose objects are the split extensions of $\mathbb{C}$ with the same fixed kernel object $X$ and whose arrows are the pairs $(g,f)$ of arrows in $\mathbb{C}$
\[\begin{tikzcd}
	X & A & B \\
	X & C & D
	\arrow["g"', from=1-2, to=2-2]
	\arrow["f", from=1-3, to=2-3]
	\arrow[equal, from=1-1,to=2-1]
	\arrow["k", from=1-1, to=1-2]
	\arrow["l"', from=2-1, to=2-2]
	\arrow["p"', shift right=1, from=1-2, to=1-3]
	\arrow["s"', shift right=1, from=1-3, to=1-2]
	\arrow["q"', shift right=1, from=2-2, to=2-3]
	\arrow["r"', shift right=1, from=2-3, to=2-2]
\end{tikzcd}\]
such that $gk=l$,$fp=qg$ and $gs=rf$.
\end{definition}

Given a group $X$, we can define a functor 
\[\begin{tikzcd}[row sep=5pt]
	{\grp^{op}} && {\mathbb{S}et} \\
	B && {\{X\rightarrow A \leftrightarrows B\}/ \sim } \\
	\\
	\\
	{B'} && {\{X\rightarrow A' \leftrightarrows B'\}/ \sim}
	\arrow["{\Act(-,X)}", from=1-1, to=1-3]
	\arrow["f"', from=2-1, to=5-1]
	\arrow[maps to, from=2-1, to=2-3]
	\arrow["{\Act(f,X)}"', from=5-3, to=2-3]
	\arrow[maps to, from=5-1, to=5-3]
\end{tikzcd}\]
where $\{X\rightarrow A \leftrightarrows B\}$ is the set of split extensions with fixed kernel object $X$ and fixed quotient object $B$; two split extensions $X\rightarrow A \leftrightarrows B$ and $X\rightarrow \overline{A} \leftrightarrows B$ are equivalent (under the equivalence relation $\sim$) if there exists an arrow $g \colon  A \rightarrow \overline{A}$ such that $gk=\overline{k}$, $gs=\overline{s}$ and $\overline{p}g=p$ 
\[\begin{tikzcd}
	X & A & B \\
	X & {\overline{A}} & B
	\arrow["{\overline{k}}"', from=2-1, to=2-2]
	\arrow["k", from=1-1, to=1-2]
	\arrow[equal, from=1-1, to=2-1]
	\arrow["g"', from=1-2, to=2-2]
	\arrow["p"', shift right=1, from=1-2, to=1-3]
	\arrow["s"', shift right=1, from=1-3, to=1-2]
	\arrow[equal, from=1-3, to=2-3]
	\arrow["{\overline{p}}"', shift right=1, from=2-2, to=2-3]
	\arrow["{\overline{s}}"', shift right=1, from=2-3, to=2-2]
\end{tikzcd}\]
(hence, thanks to the Split Short Five Lemma, $g$ is an isomorphism). Finally, $\Act(f,X)$ sends the class of a split extension $X\rightarrow A' \leftrightarrows B'$ to the class of the split extension defined via the following diagram, where the right-rightward square is a pullback:
\[\begin{tikzcd}
	X & {A' \times_{B'}B} & B \\
	X & {A'} & {B'.}
	\arrow["{k'}"', from=2-1, to=2-2]
	\arrow["{p'}"', shift right=1, from=2-2, to=2-3]
	\arrow["{s'}"', shift right=1, from=2-3, to=2-2]
	\arrow["f", from=1-3, to=2-3]
	\arrow["{\pi_B}"', shift right=1, from=1-2, to=1-3]
	\arrow["{\pi_{A'}}"', from=1-2, to=2-2]
	\arrow["{(s'f,id_B)}"', shift right=1, from=1-3, to=1-2]
	\arrow[equal, from=1-1, to=2-1]
	\arrow["{(k',0)}", from=1-1, to=1-2]
 \arrow["\lrcorner"{anchor=center, pos=0.125}, draw=none, from=1-2, to=2-3]
\end{tikzcd}\]

It is a known fact that, in $\grp$, there is a one-to-one correspondence between the set $\{X\rightarrow A \leftrightarrows B\}/ \sim$ and the set of group homomorphisms with domain $B$ and codomain the group $\Aut(X)$ of automorphisms of $X$. In fact, it turns out that the functor $\Act(-,X)$ is representable and a representing object is $\Aut(X)$. A pointed protomodular category in which the functor $\Act(-,X)$ is representable for every object $X$ is called \emph{action representable} \cite{borceux2005representability}. It can be seen that this condition is extremely strong. In fact, it emerged that many of categories of interest in classical algebra are not action representable, but do satisfy a weaker related condition we will describe now. In order to do this, it is easy to observe that, in $\grp$, the property of being an action representable category can be restated in the following way: the split extension 
\begin{equation}\label{terminal}
\begin{tikzcd}
	X & {X \rtimes \Aut(X)} & {\Aut(X),}
	\arrow["{(id_X,0)}"', from=1-1, to=1-2]
	\arrow["{\pi_{\Aut(X)}}"', shift right=1, from=1-2, to=1-3]
	\arrow["{(0,id_{\Aut(X)})}"', shift right=1, from=1-3, to=1-2]
\end{tikzcd}
\end{equation}
 corresponding to the action $id_{\Aut(X)} \colon  \Aut(X) \rightarrow \Aut(X)$, is a terminal object of $\splt{X}{\grp}$. Therefore, in $\grp$, for every object of $\splt{X}{\grp}$ there exists a unique morphism into \eqref{terminal}.\\
In light of this, the authors of \cite{artactionaccessible} have weakened the notion of action representable category in the following way:

\begin{definition}[\cite{artactionaccessible}] An object $F$ of $\splt{X}{\mathbb{C}}$ is said to be \emph{faithful} if for each object $E$ of $\splt{X}{\mathbb{C}}$ there is at most one arrow from $E$ to $F$.

\end{definition}

\begin{definition}[\cite{artactionaccessible}] Let $\cat{C}$ be a pointed protomodular category. An object in $\splt{X}{\mathbb{C}}$ is said to be \emph{accessible} if it admits a morphism into a faithful object. We say that $\cat{C}$ is \emph{action accessible} if, for every object $X$ of $\cat{C}$, every object in $\splt{X}{\mathbb{C}}$ is accessible.

\end{definition}

As mentioned above, the notion of action accessible category appears as a generalization of the one of action representable category: if there is a terminal object $ T $ of $\splt{X}{\mathbb{C}}$, this object is also faithful and each object of $\splt{X}{\mathbb{C}}$ admits a unique morphism into $T$. Examples of action accessible categories include, for instance, not necessarily unitary rings (as shown in \cite{artactionaccessible}) and all \emph{categories of interest} in the sense of \cite{orzech} (as shown in \cite{montoliorzech}).\\ 
One of the interesting properties, among other things, implied by action accessibility is that centralizers of internal equivalence relations exist and they have a simple description (see Theorem 4.1 in \cite{artactionaccessible}). Moreover, action accessibility implies (SH) (Theorem 5.4 in \cite{artactionaccessible}); hence, in an action accessible category centralizers of normal subobjects exist and they are normal. In the next part of this section we show that the category $\lgrp$ is not action accessible despite the fact that centralizers of internal equivalence relations exist and they have the same simple description as in action accessible categories (as shown at the end of previous section).\\

\begin{proposition} \label{non accessibile}

$\lgrp$ is not action accessible.

\begin{proof}

Consider the lexicographic product $\mathbb{Z} \cev{\times} \mathbb{Z}$ of the group of integers $\mathbb{Z}$ (with the usual order) with itself. The underlying set is the product, and in terms of structure the group operations are defined component-wise, while the order is defined as follows: $(a,b) \leq (c,d)$ if and only if $b<d$, or $b=d$ and $a \leq c$. We consider the following split extension
\begin{equation}\label{lessico}
    \begin{tikzcd} 
	{\mathbb{Z}} & {\mathbb{Z} \cev{\times} \mathbb{Z}} & {\mathbb{Z},}
	\arrow["{i_1}", from=1-1, to=1-2]
	\arrow["{p_2}"', shift right=1, from=1-2, to=1-3]
	\arrow["{i_2}"', shift right=1, from=1-3, to=1-2]
\end{tikzcd}
\end{equation}
where $i_1=(id_{\mathbb{Z}},0)$, $i_2=(0,id_{\mathbb{Z}})$ and $p_2$ is the projection on the second component. Now, for every $n \in \mathbb{N}_{>0}$ we can consider the morphism of lattice-ordered groups $f_n \colon  \mathbb{Z} \rightarrow \mathbb{Z}$ given by $f_n(x) \coloneqq  nx$. This morphism induces the following morphism in $\splt{\mathbb{Z}}{\lgrp}$:
\[\begin{tikzcd}
	{\mathbb{Z}} & {\mathbb{Z} \cev{\times} \mathbb{Z}} & {\mathbb{Z}} \\
	{\mathbb{Z}} & {\mathbb{Z} \cev{\times} \mathbb{Z}} & {\mathbb{Z}}
	\arrow["{i_1}", from=1-1, to=1-2]
	\arrow["{p_2}"', shift right=1, from=1-2, to=1-3]
	\arrow["{i_2}"', shift right=1, from=1-3, to=1-2]
	\arrow[equal, from=1-1, to=2-1]
	\arrow["{i_1}"', from=2-1, to=2-2]
	\arrow["{p_2}"', shift right=1, from=2-2, to=2-3]
	\arrow["{i_2}"', shift right=1, from=2-3, to=2-2]
	\arrow["{f_n}", from=1-3, to=2-3]
	\arrow["{g_n}"', from=1-2, to=2-2]
\end{tikzcd}\]
where $g_n(x,y) \coloneqq  (x,ny)$.
Thus we can deduce that ($\ref{lessico}$) is not faithful. So, if $\lgrp$ were action accessible then there should exist a faithful object
\[\begin{tikzcd}
	X & A & B
	\arrow["k"', from=1-1, to=1-2]
	\arrow["p"', shift right=1, from=1-2, to=1-3]
	\arrow["s"', shift right=1, from=1-3, to=1-2]
\end{tikzcd}\]
and a morphism
\[\begin{tikzcd}
	{\mathbb{Z}} & {\mathbb{Z} \cev{\times} \mathbb{Z}} & {\mathbb{Z}} \\
	\mathbb{Z} & A & B.
	\arrow["k"', from=2-1, to=2-2]
	\arrow["p"', shift right=1, from=2-2, to=2-3]
	\arrow["s"', shift right=1, from=2-3, to=2-2]
	\arrow["{i_1}", from=1-1, to=1-2]
	\arrow["{p_2}"', shift right=1, from=1-2, to=1-3]
	\arrow["{i_2}"', shift right=1, from=1-3, to=1-2]
	\arrow[equal, from=1-1, to=2-1]
	\arrow["g"', from=1-2, to=2-2]
	\arrow["f", shift left=1, from=1-3, to=2-3]
\end{tikzcd}\]
Then, if we consider the (regular epimorphism, monomorphism)-factorization of $(g,f)$ we get:
\begin{equation}\label{quadratone}
\begin{tikzcd}
	{\mathbb{Z}} & {\mathbb{Z} \cev{\times} \mathbb{Z}} & {\mathbb{Z}} \\
	{\mathbb{Z}} & {\Imm(g)} & {\Imm(f)} \\
	{\mathbb{Z}} & A & B.
	\arrow["k"', from=3-1, to=3-2]
	\arrow["p"', shift right=1, from=3-2, to=3-3]
	\arrow["s"', shift right=1, from=3-3, to=3-2]
	\arrow["{i_1}", from=1-1, to=1-2]
	\arrow["{p_2}"', shift right=1, from=1-2, to=1-3]
	\arrow["{i_2}"', shift right=1, from=1-3, to=1-2]
	\arrow["g"'{pos=0.2}, bend right, from=1-2, to=3-2]
	\arrow["f"{pos=0.2}, shift left=1, bend left, from=1-3, to=3-3]
	\arrow[equal, from=1-1, to=2-1]
	\arrow[equal, from=2-1, to=3-1]
	\arrow[tail, from=2-2, to=3-2]
	\arrow[two heads, from=1-2, to=2-2]
	\arrow["{\overline{k}}"', from=2-1, to=2-2]
	\arrow[two heads, from=1-3, to=2-3]
	\arrow[tail, from=2-3, to=3-3]
	\arrow["{\overline{p}}"', shift right=1, from=2-2, to=2-3]
	\arrow["{\overline{s}}"', shift right=1, from=2-3, to=2-2]
\end{tikzcd}
\end{equation}
Therefore, $\Imm(f)$ is a quotient (in $\lgrp$) of $\mathbb{Z}$. However, $\mathbb{Z}$ has only two ideals: $\{ 0 \}$ and $\mathbb{Z}$. Hence we have two possibilities: $\Imm(f) \cong \mathbb{Z}$ or $\Imm(f) \cong \{ e \}$. If $\Imm(f) \cong \mathbb{Z}$ then $f$ is injective and so the split extension 
\[\begin{tikzcd}
	{\mathbb{Z}} & {\mathbb{Z} \cev{\times} \mathbb{Z}} & {\mathbb{Z}}
	\arrow["{i_1}", from=1-1, to=1-2]
	\arrow["{p_2}"', shift right=1, from=1-2, to=1-3]
	\arrow["{i_2}"', shift right=1, from=1-3, to=1-2]
\end{tikzcd}\]
has to be a faithful object, and this is a contradiction. Alternatively, if $\Imm(f)\cong \{ e \}$, $f$ has to be the trivial morphism, and so $\Imm(g) \cong \mathbb{Z}$. Therefore, recalling that the top right-rightward square of \eqref{quadratone} is a pullback, we get a contradiction since $\mathbb{Z} \cev{\times} \mathbb{Z}$ is not isomorphic, as a lattice-ordered group, to $\mathbb{Z} \times \mathbb{Z}$.
\end{proof}

\end{proposition}

In \cite{artactionaccessible} the authors show that, in the case of the variety $\mathbb{R}ng$ of not necessarily unitary rings, given a split extension there is a procedure, based on centralizers of subobjects, to build a morphism from it into a faithful split extension. This same argument has also been extended in \cite{montoliorzech} to categories of interest in the sense of \cite{orzech}. We recall here a sketch of the proof presented in \cite{artactionaccessible}. \\Given an object $A$ of $\mathbb{R}ng$, two subobjects $X,Y \leq A$ cooperate if and only if, for every $x \in X$ and for every $y \in Y$, $$xy=0=yx.$$ Hence, it can be shown that the centralizer of $X$ in $A$ is the subobject $$Z_A(X) \coloneqq  \{ a \in A \, | \, ax=0=xa \text{ for all } x \in X \}.$$
Given an object of $\splt{X}{\mathbb{R}ng}$
\begin{equation}\label{split1}
    \begin{tikzcd}
	X & A & B,
	\arrow["k", from=1-1, to=1-2]
	\arrow["p", shift left=1, from=1-2, to=1-3]
	\arrow["s", shift left=1, from=1-3, to=1-2]
\end{tikzcd}
\end{equation}
they define $I \coloneqq  \{b \in B \, | \, s(b)k(x)=0=k(x)s(b) \text{ for all } x \in X \}$; they prove that $I$ is an ideal of $B$ and $s(I)=Z_A(k(X)) \cap s(B)$ is an ideal of $A$. Thus, they show that the split extension
\begin{equation}\label{split2}
   \begin{tikzcd}
	X & A/s(I) & B/I,
	\arrow["\overline{k}", from=1-1, to=1-2]
	\arrow["\overline{p}", shift left=1, from=1-2, to=1-3]
	\arrow["\overline{s}", shift left=1, from=1-3, to=1-2]
\end{tikzcd} 
\end{equation}
where the morphisms are induced by the universal property of the quotient, is a faithful object of $\splt{X}{\mathbb{R}ng}$ and the pair $(\pi_{s(I)},\pi_I)$, obtained by the quotient projections $\pi_{s(I)} \colon  A \rightarrow A/s(I)$ and $\pi_I \colon  B \rightarrow B/I$, is a morphism between \eqref{split1} and \eqref{split2}. Therefore, in the case of rings, there is a canonical way to construct an arrow of $\splt{X}{\mathbb{R}ng}$ into a faithful object making use of the notion of centralizer.\\
Although the category $\lgrp$ is not action accessible, it is possible to emulate the previous construction in this case. This shows that in $\lgrp$ centralizers of subobjects have a good behaviour even though the category is not action accessible.\\ Let us fix a split extension in $\lgrp$:
\[\begin{tikzcd}
	X & A & B.
	\arrow["k", from=1-1, to=1-2]
	\arrow["p", shift left=1, from=1-2, to=1-3]
	\arrow["s", shift left=1, from=1-3, to=1-2]
\end{tikzcd}\]
We want to show that the intersection between $s(B)$ and the centralizer of $k(X)$ in $A$ is an ideal of $A$. In other words, we need to prove that $\pol{k(X)} \cap S(B)$ is convex and closed under conjugation in $A$. 
\begin{itemize}
    \item Convexity: let us consider $s(b_1) \leq a \leq s(b_2)$ where $s(b_1), s(b_2) \in \pol{k(X)} \cap s(B)$ and $a \in A$. We recall that for all $a \in A$ there exist $x \in k(X)$ and $b \in B$ such that $a=k(x)s(b)$ (see Proposition \ref{semidiretto}). Then, applying $p$ to the inequalities, we obtain $b_1 \leq b \leq b_2$ and thus $$s(b_1) \leq s(b) \leq s(b_2).$$ Therefore, since $\pol{k(X)}$ is a convex subobject of $A$, we get $s(b) \in \pol{k(X)} \cap s(B)$. Hence, from $s(b_1) \leq a \leq s(b_2)$ multiplying on the right by $s(\inv{b})$, we obtain $$s(b_1\inv{b}) \leq k(x) \leq s(b_2\inv{b}).$$ So, since $s(b_1\inv{b}),s(b_2\inv{b}) \in \pol{k(X)} \cap s(B)$ (because $\pol{k(X)} \cap s(B)$ is a subalgebra of $A$ and $s(b), s(b_1), s(b_2) \in \pol{k(X)} \cap s(B)$), we get $k(x) \in \pol{k(X)}$ ($\pol{k(X)}$ is a convex subalgebra). Therefore, $k(x)=e$ and then we obtain $a=k(x)s(b)=s(b) \in \pol{k(X)} \cap s(B)$.
    \item Closedness under conjugation: let us consider $s(c) \in \pol{k(X)} \cap s(B)$ and $a=k(x)s(b) \in A$. Then, we have $$as(c) \inv{a}=k(x)s(b)s(c)s(\inv{b})\inv{k({x})}=k(x)s(d)\inv{k(x)}$$ where $s(d)=s(b)s(c)s(\inv{b}) \in \pol{k(X)} \cap s(B)$ since $\pol{k(X)}$ is closed under conjugation and, clearly, $s(d) \in s(B)$. Therefore, we get $$as(c) \inv{a}=k(x)s(d)\inv{k(x)}=s(d)$$ since $s(d) \in \pol{k(X)}$ and the elements of $\pol{k(X)}$ commute with the ones of $k(X)$.
\end{itemize}

\section{Fiber-wise Algebraic Cartesian Closedness} \label{Fiber-wise Algebraic Cartesian Closedness}

In this section we deal with a stronger version of the notion of algebraically cartesian closed category. We propose an equivalent version of the definition presented in \cite{acc}:

\begin{definition}[\cite{acc}] A category $\mathbb{C}$ is \emph{fiber-wise algebraically cartesian closed} if for every split epimorphism 
\[\begin{tikzcd}
	A & B
	\arrow["p", shift left=1, from=1-1, to=1-2]
	\arrow["s", shift left=1, from=1-2, to=1-1]
\end{tikzcd}\]
the change-of-base functor $$p^* \colon  \pt{B}{\mathbb{C}} \rightarrow \pt{A}{\mathbb{C}}$$ has a right adjoint.

\end{definition}

As established in \cite{acc}, it is not difficult to see that this condition holds for a category $\cat{C}$ if and only if every category of points over $\cat{C}$ is algebraically cartesian closed. First of all, we observe that the category $\pt{A \leftrightarrows B}{\pt{B}{\mathbb{C}}}$ is isomorphic to $\pt{A}{\mathbb{C}}$: an object of $\pt{A \leftrightarrows B}{\pt{B}{\mathbb{C}}}$ can be seen as a diagram of type 
\[\begin{tikzcd}
	C & B \\
	A & B,
	\arrow["q"', shift right=1, from=1-1, to=1-2]
	\arrow["r"', shift right=1, from=1-2, to=1-1]
	\arrow["p"', shift right=1, from=2-1, to=2-2]
	\arrow["s"', shift right=1, from=2-2, to=2-1]
	\arrow[equal, from=1-2, to=2-2]
	\arrow["h"', shift right=1, from=1-1, to=2-1]
	\arrow["t"', shift right=1, from=2-1, to=1-1]
\end{tikzcd}\]
where $ps=id_B$, $qr=id_B$, $ht=id_A$, $ph=q$, and $ts=r$; therefore, $q$ and $r$ are uniquely determined by $h$ and $t$. Hence, each category of points is algebraically cartesian closed if and only if the functor
\[\begin{tikzcd}
	{\tau^* \colon  \pt{B=B}{\pt{B}{\mathbb{C}}}} & {\pt{A \leftrightarrows B}{\pt{B}{\mathbb{C}}}}
	\arrow[from=1-1, to=1-2]
\end{tikzcd}\]
has a right adjoint; thanks to the isomorphism shown above between $\pt{A}{\mathbb{C}}$ and $\pt{A \leftrightarrows B}{\pt{B}{\mathbb{C}}}$ and recalling that $\tau=p$, we get that $\tau^*$ has a right adjoint if and only if $p^*$ has a right adjoint. Therefore, $\mathbb{C}$ is fiber-wise algebraically cartesian closed if and only if every category of points over $\cat{C}$ is algebraically cartesian closed.\\
Our aim is to show, thanks to the previous observations, that $\lgrp$ is fiber-wise algebraically cartesian closed. In order to do this we will prove that, in every category of points over $\lgrp$, subobjects have centralizers. As a preliminary remark, we recall that an arrow $(A \leftrightarrows B) \overset{f}{\rightarrow} (C \leftrightarrows B)$ of $\pt{B}{\mathbb{C}}$ is a monomorphism if and only if $f \colon A \rightarrow C$ is a monomorphism of $\mathbb{C}$.\\

We are ready to show the existence of centralizers in every category of points over $\lgrp$ and to provide an explicit description of them. 

\begin{definition} \label{chiusura}

Let $X$ be an object of $\lgrp$ and $B$ a subalgebra of $X$. A subalgebra $L$ of $X$ is \emph{closed under the action of} $B$ if $$\jnoc{b}{l} \in L \textit{ and } (l_1b_1 \lor l_2b_2) \inv{(b_1 \lor b_2)} \in L$$ for every $l,l_1,l_2 \in L$ and $b,b_1,b_2 \in B$.

\end{definition}

\begin{proposition}\label{punti e centralizzanti}

Let $B$ be an object of $\lgrp$. In the category $\pt{B}{\lgrp}$ subobjects have centralizers.

\begin{proof}
Let us consider an object $(A,p,s)$ of $Pt_B \lgrp$, i.e.\ a diagram of the form
\[\begin{tikzcd}
	K & A & B
	\arrow["p"', shift right=1, from=1-2, to=1-3]
	\arrow["s"', shift right=1, from=1-3, to=1-2]
	\arrow["k", from=1-1, to=1-2]
\end{tikzcd}\]
where $ps=id_B$ and $k=\kerr(p)$. Given two subobjects $(X,p_{|X},r_X)$ and $(Y,p_{|Y},r_Y)$ of $(A,p,s)$ and the product between them in $\pt{B}{\lgrp}$, we need to describe the arrows $i_X \colon X \rightarrow X \times_B Y$ and $i_Y \colon Y \rightarrow X \times_B Y$ induced by the universal property (we recall that the product in $\pt{B}{\lgrp}$ is given by the pullback of $p_{|X}$ along $p_{|Y}$ in $\lgrp$). Then, if we consider the following diagram
\[\begin{tikzcd}
	X & B \\
	& {X \times_BY} & Y \\
	& X & B,
	\arrow["{p_{|Y}}", from=2-3, to=3-3]
	\arrow["{p_{|X}}"', from=3-2, to=3-3]
	\arrow["{\pi_X}"', from=2-2, to=3-2]
	\arrow["{\pi_Y}", from=2-2, to=2-3]
	\arrow["\lrcorner"{anchor=center, pos=0.125}, draw=none, from=2-2, to=3-3]
	\arrow["{id_X}"', bend right, from=1-1, to=3-2]
	\arrow["{i_X}", from=1-1, to=2-2]
	\arrow["{p_{|X}}", from=1-1, to=1-2]
	\arrow["s", bend left, from=1-2, to=2-3]
\end{tikzcd}\]
we get $i_X(x)=(x,sp(x))$; in a similar way $i_Y \colon  Y \rightarrow X \times_B Y$ is given by $i_Y(y)=(sp(y),y)$.\\
We know that $(x,y) \in X \times_B Y$ if and only if $p(x)=p(y)$. Moreover, given $(x,y) \in X \times_B Y$, we have $(x,y)=(x,sp(x))\inv{(sp(y),sp(x))}(sp(y),y)$ where $(x,sp(x))\in X \times_B Y$, $(sp(y),y)\in X \times_B Y$ and, since $psp(x)=p(x)=p(y)=psp(y)$, we get $(sp(y),sp(x))=(sp(x),sp(x))=(sp(y),sp(y)) \in X \times_B Y$.\\ Hence, if there exists a cooperator $\varphi \colon  X \times_B Y \rightarrow A$, then $$\varphi(x,y)=\varphi(x,sp(x))\inv{\varphi(sp(y),sp(x))}\varphi(sp(y),y)=x \inv{sp(x)} y= x \inv{sp(y)} y,$$
since $$\varphi(x,sp(x))\inv{\varphi(sp(y),sp(x))}\varphi(sp(y),y)=\varphi(i_X(x))\inv{\varphi(i_X(sp(x)))}\varphi(i_Y(y))$$
\[\begin{tikzcd}
	X && Y \\
	& {X \times_BY} \\
	& A.
	\arrow["{i_Y}"', from=1-3, to=2-2]
	\arrow["{i_X}", from=1-1, to=2-2]
	\arrow["\varphi"', from=2-2, to=3-2]
	\arrow[bend left, tail, from=1-3, to=3-2]
	\arrow[bend right, tail, from=1-1, to=3-2]
\end{tikzcd}\]
It is well known that the category of points over $B$, in the case of groups, is equivalent to the category of groups with an action of $B$. Therefore, thanks to Proposition \ref{semidiretto}, we can extend this to the case of lattice-ordered groups and obtain that $A$ is isomorphic to $K \rtimes B$ with join operation defined by $$(k_1,b_1)\lor(k_2,b_2)=((k_1b_1 \lor k_2b_2)\inv{(b_1 \lor b_2)},b_1 \lor b_2);$$ more specifically, the split extension $(A, p, s)$ is isomorphic to the split extension
\[\begin{tikzcd}
	 {K \rtimes B} & B.
	\arrow["{p_B}"', shift right=1, from=1-1, to=1-2]
	\arrow["{i_B}"', shift right=1, from=1-2, to=1-1]
\end{tikzcd}\]
Hence, a subobject $(X,q,r)$ of $(A,p,s)$ in $\pt{B}{\lgrp}$ can be seen, modulo isomorphisms, as a subalgebra $X \leq K \rtimes B$ in $\lgrp$ such that, referring to the diagram
\[\begin{tikzcd}
	X \\
	{K \rtimes B} & B,
	\arrow[hook, from=1-1, to=2-1]
	\arrow["{p_B}", shift left=1, from=2-1, to=2-2]
	\arrow["{i_B}", shift left=1, from=2-2, to=2-1]
	\arrow["q", shift left=3, from=1-1, to=2-2]
	\arrow["r", shift right=1, from=2-2, to=1-1]
\end{tikzcd}\]
$q$ is the restriction of $p_B$ to $X$ and, for every $b \in B$, $r(b)=(e,b)$ (in particular $\{e \} \times B \leq X$).\\
Given a subobject $(X,q,r)$ of $(K \rtimes B,p_B,i_B)$ we define $$\overline{X} \coloneqq  \{ \overline{x} \in K \, | \, \text{there exists } b \in B \text{ s.t. } (\overline{x},b) \in X \}.$$ 
We show that $X=\overline{X} \times B$ as sets. Clearly $X \subseteq \overline{X} \times B$. Conversely, fix an element $(k,b) \in \overline{X} \times B$. Then $k \in \overline{X}$, and so there exists $b_1 \in B$ such that $(k,b_1) \in X$; but $(e,b),(e,b_1) \in X$, therefore $(k,b)=(k,b_1)\inv{(e,b_1)}(e,b) \in X$. In general, we have a one-to-one correspondence between the subobjects of $(A,p,s)$ and the subalgebras of $K$ closed under the action of $B$.\\
Therefore, given two subobjects $(X,p_{|X},r_X)$ and $(Y,p_{|Y},r_Y)$ of $(A,p,s)$ we can suppose $X=\overline{X} \times B$ and $Y=\overline{Y} \times B$. Thus, $$X \times_B Y= \{ ((\overline{x},b_1),(\overline{y},b_2)) \in X \times Y \, | \, b_1=b_2 \}$$ and $\varphi \colon  X \times_B Y \rightarrow K \rtimes B$ is such that $$\varphi((\overline{x},b),(\overline{y},b))=(\overline{x},b)\inv{(e,b)}(\overline{y},b)=(\overline{x} \, \overline{y},b).$$
We start by showing that $$\varphi \text{ is a group homomorphism if and only if }xy=yx \text{ for all }x \in \overline{X},y \in \overline{Y}.$$ Given an element $((x,b),(y,b)),((z,c),(w,c)) \in ( \overline{X} \times B) \times_B (\overline{Y} \times B)$ one has 
\begin{align*}
    ((x,b),(y,b))((z,c),(w,c))&=((x,b)(z,c),(y,b)(w,c))\\
    &=((x \jnoc{b}{z},bc),(y \jnoc{b}{w},bc)).
\end{align*}
Hence $$\varphi((x,b),(y,b))\varphi((z,c),(w,c))=(xy,b)(zw,c)=(xy\jnoc{b}{zw},bc)$$ and, since $((x,b),(y,b))((z,c),(w,c))=((x \jnoc{b}{z},bc),(y \jnoc{b}{w},bc))$, we get $$\varphi((x \jnoc{b}{z},bc),(y \jnoc{b}{w},bc))=(x \jnoc{b}{z}y \jnoc{b}{w},bc).$$ So, $\varphi$ is a group homomorphism if and only if $\jnoc{b}{z}y=y\jnoc{b}{z}$ for all $z \in \overline{X},\text{ } y \in \overline{Y}$ and $b \in B$. Then, setting $b=e$, we obtain $zy=yz$ for all $z \in \overline{X}$ and $y \in \overline{Y}$. Moreover, since $\overline{X}$ is closed under the action of $B$ and the conjugation is a bijection, we get that every element of $\overline{X}$ can be seen as $\jnoc{b}{z}$, for appropriate $b \in B$ and $z \in \overline{X}$; thus, if $zy=yz$ for all $z \in \overline{X}$ and $y \in \overline{Y}$ then $\jnoc{b}{z}y=y\jnoc{b}{z}$ for all $z \in \overline{X}, \text{ } y \in \overline{Y}$ and $b \in B$.
Now, let us deal with the order structure. We know that $\varphi$ is a morphism of lattice-ordered groups if and only if $\varphi$ is a group homomorphism and for all $((x,b),(y,b)) \in (\overline{X} \times B) \times_B (\overline{Y} \times B)$
\begin{equation} \label{equazione}
\varphi(((x,b),(y,b)) \lor ((e,e),(e,e)))=\varphi(((x,b),(y,b))) \lor (e,e).
\end{equation}
Observing that
\begin{align*}
((x,b),(y,b)) \lor ((e,e),(e,e))&=((x,b) \lor (e,e),(y,b) \lor (e,e))\\
&=(((xb \lor e) \inv{(b \lor e)},b \lor e),((yb \lor e) \inv{(b \lor e)}, b \lor e))\\
&=(((xb)^+\inv{(b^+)},b^+),((yb)^+\inv{(b^+)},b^+)),
\end{align*}
we have $$\varphi(((x,b),(y,b)) \lor ((e,e),(e,e)))=((xb)^+\inv{(b^+)}(yb)^+\inv{(b^+)},b^+).$$ Considering the right term of \eqref{equazione}, we obtain $$\varphi(((x,b),(y,b))) \lor (e,e)=(xy,b) \lor (e,e)=((xyb)^+ \inv{(b^+)},b^+).$$
We prove that $$\varphi \text{ is a lattice-ordered group morphism if and only if } \overline{X} \perp \overline{Y}.$$
If $\varphi$ is a lattice-ordered group morphism, then $$(xb)^+\inv{(b^+)}(yb)^+\inv{(b^+)}=(xyb)^+\inv{(b^+)}$$ for each $x \in \overline{X}$, $y \in \overline{Y}$ and $b \in B$; therefore, setting $b=e$, we get $$(xy)^+=x^+y^+ \text{ for every } x \in \overline{X} \text{ and } y \in \overline{Y},$$ and so, thanks to Corollary \ref{magie tra ortogonali}, $\overline{X}$ and $\overline{Y}$ are orthogonal.\\
Conversely, let us suppose $\overline{X} \perp \overline{Y}$. We want to show that, for all $x \in \overline{X},y \in \overline{Y} \text{ and }b \in B$,  $$(xb)^+\inv{(b^+)}(yb)^+\inv{(b^+)}=(xyb)^+\inv{(b^+)}$$ or, equivalently, that $(xb \lor e)\inv{(b \lor e)}(yb \lor e)=xyb \lor e$. We start with the term on the left: 
\begin{align*}
(xb \lor e)\inv{(b \lor e)}(yb \lor e)&=(xb( \inv{b} \land e) \lor (\inv{b} \land e))(yb \lor e)\\
&=xb(\inv{b} \land e)(yb \lor e) \lor (\inv{b} \land e)(yb \lor e)\\
&=xb(\inv{b} \land e)yb \lor xb(\inv{b} \land e) \lor (\inv{b} \land e)yb \lor (\inv{b} \land e).
\end{align*}
Now, we know that there exists an element $y_1 \in \overline{Y}$ such that $yb=by_1$ (since $\overline{Y}$ is closed under the action of $B$ and the conjugation is an automorphism of $\overline{Y}$), hence the last term is equal to 
\begin{align*}
xb(\inv{b} \land e)by_1 &\lor xb(\inv{b} \land e) \lor (\inv{b} \land e)by_1 \lor (\inv{b} \land e)\\
&=x(b \land b^2)y_1 \lor x(b \land e) \lor (b \land e)y_1 \lor (\inv{b} \land e).
\end{align*}
Moreover, we observe that $(b \land e)y_1=y_2(b \land e)$ for an appropriate element $y_2 \in \overline{Y}$, thus 
\begin{align*}
x(b \land e) \lor (b \land e) y_1&=x(b \land e) \lor y_2(b \land e)=(x \lor y_2)(b \land e)\\
&=(xy_2 \lor e)(b \land e)=xy_2(b \land e) \lor (b \land e)\\
&=x(b \land e)y_1 \lor (b \land e)
\end{align*}
(by Corollary \ref{magie tra ortogonali} we know that $x \lor y_2=xy_2 \lor e$, since $\overline{X} \perp \overline{Y}$). Then, one has \begin{align*}
x(b \land b^2)y_1 &\lor x(b \land e) \lor (b \land e)y_1 \lor (\inv{b} \land e)\\
&=x(b \land b^2)y_1 \lor x(b \land e)y_1 \lor (b \land e) \lor (\inv{b} \land e).
\end{align*}
We recall that $(b \land e) \lor (\inv{b} \land e)=(b \lor \inv{b}) \land e=|b| \land e=e$, so we finally get 
\begin{align*}
    x(b \land b^2)y_1 \lor x(b \land e)y_1 \lor (b \land e) \lor (\inv{b} \land e)&=x(b \land b^2)y_1 \lor x(b \land e)y_1 \lor e\\
&=x[(b \land b^2) \lor (b \land e)]y_1 \lor e\\
&=xby_1 \lor e=xyb \lor e,
\end{align*} observing that $(b \land b^2) \lor (b \land e)=b(b \land e) \lor e(b \land e)=(b \lor e )(b \land e)=b^+b^-=b.$\\
To conclude, we have to prove that, for every subalgebra $\overline{X} \leq K$ closed under the action of $B$, then $\pol{\overline{X}} \leq K$ is closed under the action of $B$ (and so the centralizer of $X=\overline{X} \times B$ is $\pol{\overline{X}} \times B$ endowed with the structure induced by the semi-direct product). Fix an element $y \in \pol{\overline{X}}$ and an element $b \in B$; we show that $\jnoc{b}{y} \in \pol{\overline{X}}$. We know that $$\displaylines{\jnoc{b}{y} \in \pol{\overline{X}} \iff |\jnoc{b}{y}| \land |x|=e \text{ for all } x \in \overline{X} \iff\cr |y| \land |\conj{b}{x}|=e \text{ for all } x \in \overline{X} \iff |y| \land |x|=e \text{ for all } x \in \overline{X};}$$ observe that the last assertion holds since $y \in \pol{\overline{X}}$. We recall that $\pol{\overline{X}}$ is a convex subalgebra of $K$. For every $y_1,y_2 \in \pol{\overline{X}}$ and $b_1,b_2 \in B$ we have $y_1b_1 \leq (y_1 \lor y_2)(b_1 \lor b_2)$ and $y_2b_2 \leq (y_1 \lor y_2)(b_1 \lor b_2)$; we also observe that $(y_1 \land y_2)b_1 \leq y_1b_1$ and $(y_1 \land y_2)b_2 \leq y_2b_2$. So one has $$y_1b_1 \lor y_2b_2 \leq (y_1 \lor y_2)(b_1 \lor b_2)$$ and $$(y_1 \land y_2)(b_1 \lor b_2)=(y_1 \land y_2)b_1 \lor (y_1 \land y_2)b_2 \leq y_1b_1 \lor y_2b_2.$$ Therefore, for all $y_1,y_2 \in \pol{\overline{X}}$ and $b_1,b_2 \in B$, we obtain $$y_1 \land y_2 \leq (y_1b_1 \lor y_2b_2) \inv{(b_1 \lor b_2)}\leq y_1 \lor y_2;$$ then, since $\pol{\overline{X}}$ is convex in $K$, we get $(y_1b_1 \lor y_2b_2)\inv{(b_1 \lor b_2)} \in \pol{\overline{X}}$ for all $y_1,y_2 \in \pol{\overline{X}}$ and $b_1,b_2 \in B$.
\end{proof}
\end{proposition}

\begin{corollary}

For every object $B$ of $\lgrp$, in the category $\pt{B}{\lgrp}$ subobjects have centralizers, therefore $\pt{B}{\lgrp}$ is algebraically cartesian closed. Hence, the category $\lgrp$ is fiber-wise algebraically cartesian closed.

\end{corollary}

\section{Normality of the Higgins Commutator} \label{Normality of the Higgins Commutator}

The aim of this section is to propose a further study regarding the properties of commutators in the category of lattice-ordered groups.\\
We recall a first notion of categorical commutator strongly linked to the concept of cooperation.

\begin{definition}[\cite{huq}, \cite{cooperare}] Let $\mathbb{C}$ be a unital category. For a pair of subobjects $a \colon  A \mono X$ and $b \colon  B \mono X$ of an object $X$ in $\mathbb{C}$,
the \emph{Huq commutator} is the smallest normal subobject $[A,B]_X \mono X$ such that the
images of a and b cooperate in the quotient $X/[A,B]_X$.

\end{definition}

Next we recall here the notion of \emph{Higgins commutator}. In a pointed category $\mathbb{C}$ with binary products and coproducts, for every pair of objects $H$ and $K$, we have the following canonical arrows:

\[\begin{tikzcd}
	H & {H \times K} & K \\
	H & {H + K} & K;
	\arrow["{(id_H,0)}", from=1-1, to=1-2]
	\arrow["{(0,id_K)}"', from=1-3, to=1-2]
	\arrow["{[0,id_K]}", from=2-2, to=2-3]
	\arrow["{[id_H,0]}"', from=2-2, to=2-1]
\end{tikzcd}\]
combining them we get a canonical arrow $$\Sigma=\big(\begin{smallmatrix} id_H & 0\\ 0 & id_K \end{smallmatrix}\big) \colon  H + K \rightarrow H \times K.$$
In other words, $\Sigma$ is the unique morphism making either the diagram
\[\begin{tikzcd}
	H & {H + K} & K \\
	& {H \times K}
	\arrow["{\iota_K}", from=1-1, to=1-2]
	\arrow["{\iota_H}"', from=1-3, to=1-2]
	\arrow["{(id_H,0)}"', from=1-1, to=2-2]
	\arrow["{(0,id_K)}", from=1-3, to=2-2]
	\arrow["\Sigma", from=1-2, to=2-2]
\end{tikzcd}\]
or equivalently the diagram
\[\begin{tikzcd}
	& {H + K} \\
	H & {H \times K} & K
	\arrow["{\pi_H}", from=2-2, to=2-1]
	\arrow["{\pi_K}"', from=2-2, to=2-3]
	\arrow["{[id_H,0]}"', from=1-2, to=2-1]
	\arrow["{[0,id_K]}", from=1-2, to=2-3]
	\arrow["\Sigma", from=1-2, to=2-2]
\end{tikzcd}\]
commute. For instance, in the case of the variety of groups, the morphism $\Sigma$ associates to each word $h_1k_1h_2k_2 \dots h_nk_n$, where $h_i \in H$ and $k_i \in K$ for $i=1, \dots, n$, the pair $(h_1h_2 \dots h_n,k_1k_2 \dots k_n) \in H \times K$. It is easy to see that a category with binary products and coproducts is unital if and only if, for every pair of objects $H$ and $K$, $\Sigma$ is a strong epimorphism (see e.g.\ \cite{semi}). Hence, again in the case of groups, the kernel of $\Sigma$, denoted by $H \diamond K$ and called the \emph{cosmash product} of $H$ and $K$, can be described as the subgroup of $H + K$ generated by  the elements of the form $hk\inv{h}\inv{k}$ with $h \in H$ and $k \in K$. In the light of the above, we are ready to recall the following:

\begin{definition}[\cite{artcommutatorienormalita}] Let $\mathbb{C}$ be a semi-abelian category. Given a pair of subobjects $a \colon  A \mono X$ and $b \colon  B \mono X$ of an object $X$ in $\mathbb{C}$, the \emph{Higgins commutator} of A and B is the subobject $[A,B] \mono X$
constructed, via the (regular epimorphism, monomorphism)-factorization, as in diagram

\[\begin{tikzcd}
	{A \diamond B} & {A+B} \\
	{[A,B]} & X,
	\arrow["{[a,b]}", from=1-2, to=2-2]
	\arrow["{k_{A,B}}", from=1-1, to=1-2]
	\arrow[two heads, from=1-1, to=2-1]
	\arrow[tail, from=2-1, to=2-2]
\end{tikzcd}\]
where $k_{A,B}$ is the kernel of $\Sigma=\big(\begin{smallmatrix} id_A & 0\\ 0 & id_B \end{smallmatrix}\big) \colon  A+B \rightarrow A \times B$. 

\end{definition}

In general, the Higgins commutator of two normal subobjects is not normal. Therefore, it makes sense to mention the following definition:

\begin{definition}[\cite{NH}] A semi-abelian category $\mathbb{C}$ satisfies the condition of \emph{normality of Higgins commutators (NH)} when, for every pair of normal subobjects $H \mono X$, $K \mono X$ where $X$ is an object of $\mathbb{C}$, the Higgins commutator $[H,K] \mono X$ is a normal subobject of $X$.

\end{definition}

We have everything we need to prove that that the category of lattice-ordered groups satisfies (NH).

\begin{lemma} \label{NH}

Let $H,K \leq X$ be two convex subalgebras of $X$ in $\lgrp$. Then $H$ and $K$ cooperate if and only if $H \cap K=\{e\}$.

    \begin{proof}
    
    ($\Rightarrow$) Trivial since $H \perp K$ (thanks to Proposition \ref{central}).\\
    ($\Leftarrow$) We want to show that $H \perp K$: let us consider two elements $h \in H$ and $k \in K$; then $e \leq |h| \land |k| \leq |h|$ and $e \leq |h| \land |k| \leq |k|$. Therefore, since $H$ and $K$ are convex, we have $|h| \land |k| \in H \cap K= \{e\}$.
    \end{proof}

\end{lemma}

\begin{notation}

We will write $[H,K]$ for the Higgins commutator of $H \mono X$ and $K \mono X$, and $[H,K]_Y$ for the Huq commutator of $H \mono X$ and $K \mono X$ in the subobject $Y$ of  $X$, where $H$ and $K$ are subobjects of $Y$. 

\end{notation}

\begin{proposition}

Let $X$ be an object of $\lgrp$ and $H,K$ ideals of $X$. Then, $[H,K]_X=H \cap K$.

\begin{proof}

Let us prove the inclusion $H \cap K \subseteq [H,K]_X$. Consider the following diagram:  
\[\begin{tikzcd}
	H && K \\
	& {H \times K} \\
	& {X/[H,K]_X}
	\arrow["\varphi", from=2-2, to=3-2]
	\arrow["{q_{|H}}"', bend right, from=1-1, to=3-2]
	\arrow["{q_{|K}}", bend left, from=1-3, to=3-2]
	\arrow["{i_H}", from=1-1, to=2-2]
	\arrow["{i_K}"', from=1-3, to=2-2]
\end{tikzcd}\]
where $q \colon X \twoheadrightarrow X/[H,K]_X$ is the canonical projection. Then, by Lemma \ref{NH}, we know that $q(H) \cap q(K)= \{e \}$. So, since $q(H \cap K) \subseteq q(H) \cap q(K)= \{e\}$ we get $H \cap K \subseteq [H,K]_X$. The other inclusion holds in every semi-abelian category (see Theorem 3.9 in \cite{everaert2012relative}).
\end{proof}

\end{proposition}

\begin{proposition}

The category $\lgrp$ satisfies (NH).

\begin{proof}

Thanks to Theorem 2.8 of \cite{NH}, it suffices to prove that, given an ideal $H$ of $X$ and an ideal $K$ of $Y$ such that $H,K \leq Y$, then $[H,K]_X=[H,K]_Y$. Thus the statement follows from the previous proposition, since $[H,K]_X=[H,K]_Y=H \cap K$.
\end{proof}

\end{proposition}

\section{Algebraic Coherence for $\lab$} \label{Algebraic Coherence}

In the last part of the paper we focus on the notion of  \emph{algebraically coherent category}. This concept has an important algebraic meaning: an algebraically coherent category satisfies a large set of properties related to the good behaviour of commutators (such as, for example, strong protomodularity); moreover, in the case of a variety of universal algebras, the property of being fiber-wise algebraically cartesian closed is implied by algebraic coherence.

\begin{definition}[\cite{artalgebricallycoherent}] A category $\mathbb{C}$ with finite limits is \emph{algebraically coherent} if, for every morphism $f \colon  X \rightarrow Y$ in $\mathbb{C}$, the change-of-base functor $$f^* \colon  \pt{Y}{\mathbb{C}} \rightarrow \pt{X}{\mathbb{C}}$$ is coherent: a functor between categories with finite limits is \emph{coherent} if it preserves finite limits and jointly extremally epimorphic pairs.
\end{definition}

Since in the semi-abelian case the split extensions with fixed splitting can be totally described in terms of semi-direct products, the authors in \cite{artalgebricallycoherent} proved the following result:

\begin{proposition}[\cite{artalgebricallycoherent}, Theorem 3.21] Suppose $\mathbb{C}$ is a semi-abelian category. The following are equivalent:
\begin{itemize}
    \item $\mathbb{C}$ is algebraically coherent;
    \item given $K \mono X$ and $H \mono X$ in $\mathbb{C}$, any action $\xi \colon  B \flat X \rightarrow X$ which restricts to $K$ and $H$ also restricts to $K \lor H$.
\end{itemize}
\end{proposition}

Let us now try to understand how this result can be interpreted in the category of lattice-ordered groups. We know (see \cite{bourn1998protomodularity}) that, in the semi-abelian case, for each internal action $\xi \colon  B \flat X \rightarrow X$ there exists a unique (up to isomorphism) split extension 

\[\begin{tikzcd}
	X & A & B
	\arrow["p"', shift right=1, from=1-2, to=1-3]
	\arrow["s"', shift right=1, from=1-3, to=1-2]
	\arrow["k", from=1-1, to=1-2]
\end{tikzcd}\]

\noindent making, in the following diagram, the right-rightward square, the right-leftward square and the left square commute:

\[\begin{tikzcd}
	{B \flat X} & {X+B} & B \\
	X & A & B.
	\arrow["{k_0}", from=1-1, to=1-2]
	\arrow["{[0,id_B]}"', shift right=1, from=1-2, to=1-3]
	\arrow["{\iota_B}"', shift right=1, from=1-3, to=1-2]
	\arrow["\xi"', from=1-1, to=2-1]
	\arrow["k"', from=2-1, to=2-2]
	\arrow["{[k,s]}"', from=1-2, to=2-2]
	\arrow["p"', shift right=1, from=2-2, to=2-3]
	\arrow["s"', shift right=1, from=2-3, to=2-2]
	\arrow[shift left=1, no head, from=1-3, to=2-3]
	\arrow[no head, from=1-3, to=2-3]
\end{tikzcd}\]

Therefore, in $\lgrp$, $\xi$ restricts to a subalgebra $L \leq X$ if and only if $L$ is closed under the corresponding action of $B$ (in the sense of Definition \ref{chiusura}).\\ 

In the next proposition we will deal with the category of \emph{lattice-ordered abelian groups}. A lattice-ordered abelian group is a lattice-ordered group in which the group operation is commutative. The category $\lab$ is the full subcategory of $\lgrp$ whose objects are lattice-ordered abelian groups.

\begin{proposition}

$\lab$ is algebraically coherent.

\begin{proof}

Let us consider an object $(A,p,s)$ of $\pt{B}{\lab}$ i.e.\ a diagram of the form
\[\begin{tikzcd}
	X & A & B
	\arrow["p"', shift right=1, from=1-2, to=1-3]
	\arrow["s"', shift right=1, from=1-3, to=1-2]
	\arrow["k", from=1-1, to=1-2]
\end{tikzcd}\]
where $ps=id_B$ and $k=\kerr(p)$. For simplicity, let us suppose that $k$ is the inclusion of $X \leq A$ and $s$ is the inclusion of $B \leq A$. Given two subalgebras $K,H \leq X$ closed under the action of $B$, we want to show that also $K \lor H$ is closed under the action of $B$.\\ 
First of all, let us observe that, given a subalgebra $L \leq X$, the following equality holds for every $l_1,l_2 \in L$ and $b_1,b_2 \in B$: $$(l_1b_1 \lor l_2b_2)\inv{(b_1 \lor b_2)}=l_2(\inv{l_2}l_1b_1\inv{b_2} \lor e)\inv{(b_1\inv{b_2}\lor e)}.$$ Therefore, $L$ is closed under the action of $B$ if and only if, for all $l \in L$ and $b \in B$, $$(lb \lor e)\inv{(b \lor e)} \text{ belongs to } L.$$
We recall, as proved in Section \ref{Preliminaries}, that in a lattice-ordered abelian group $A$ the equation $$xy=(x \lor y)(x \land y)$$ holds for all $x,y \in A$.
Finally, it is easy to see that every element of $K \lor H$ can be written as $$\bigvee_{i \in I} \bigwedge_{j \in J}k_{i,j}h_{i,j}$$ where $I,J$ are finite sets of indices and $k_{i,j} \in K,h_{i,j} \in H$ for all $i \in I, j \in J$. This statement can be proved by iteratively applying the distributive properties of the lattice operations, the distributivity property of the group product over the lattice operations, and the commutative property of the group product. Therefore, given an element $b \in B$, one has 
\begin{align*}
\Bigg( \Bigg(\bigvee_{i \in I} \bigwedge_{j \in J}k_{i,j}h_{i,j} \Bigg)b \lor e \Bigg)&=\Bigg( \bigvee_{i \in I} \bigwedge_{j \in J}k_{i,j}h_{i,j}b \Bigg) \lor e\\
&=\bigvee_{i \in I} \Bigg( \bigwedge_{j \in J}k_{i,j}h_{i,j}b \lor e \Bigg)
=\bigvee_{i \in I} \bigwedge_{j \in J}(k_{i,j}h_{i,j}b \lor e)
\end{align*}
where the first equality holds thanks to the distributivity of the group operation over the lattice operations, the second thanks to the idempotence of the join, and the third thanks to the distributivity of the join over the meet. Therefore $$\Bigg( \Bigg(\bigvee_{i \in I} \bigwedge_{j \in J}k_{i,j}h_{i,j} \Bigg)b \lor e \Bigg)\inv{(b \lor e)}=\bigvee_{i \in I} \bigwedge_{j \in J}(k_{i,j}h_{i,j}b \lor e)\inv{(b \lor e)}$$ and hence, in order to prove that if $K$ and $H$ are closed under the action of $B$ then also $K \lor H$ is closed under the action of $B$, it suffices to prove that $$(khb \lor e) \inv{(b \lor e)} \in K \lor H,$$
for every $k \in K$, $h \in H$ and $b \in B$.
We need to take care of an intermediate step: we want to show that, for all $k \in K$, $h \in H$ and $b \in B$, $$((k \lor h)b \lor \inv{(k \land h)})\inv{(b \lor e)} \text{ belongs to } K \lor H.$$ To do this we observe that 
\begin{align*}
((k \lor h)b \lor \inv{(k \land h)})\inv{(b \lor e)}&=(kb \lor hb \lor \inv{k} \lor \inv{h})\inv{(b \lor e)}\\
&=(kb \lor \inv{k})\inv{(b \lor e)} \lor (hb \lor \inv{h})\inv{(b \lor e)}\\
&=\inv{k}(k^2 b \lor e) \inv{(b \lor e)} \lor \inv{h} (h^2 b \lor e) \inv{(b \lor e)}.
\end{align*}
So, since $K$ is closed under the action of $B$, we have $(k^2 b \lor e) \inv{(b \lor e)} \in K$ and then we obtain 
$\inv{k}(k^2 b \lor e) \inv{(b \lor e)} \in K$; similarly $\inv{h} (h^2 b \lor e) \inv{(b \lor e)} \in H$. Therefore, taking the join of these two terms, we get $((k \lor h)b \lor \inv{(k \land h)})\inv{(b \lor e)} \in K \lor H.$ To conclude, since $k \land h \in K \lor H$, the product $$(k \land h)((k \lor h)b \lor \inv{(k \land h)})\inv{(b \lor e)} \text{ belongs to } K \lor H;$$ then, applying the distributivity property of the group product over the lattice join and recalling that $(k \land h)(k \lor h)=kh$, we deduce 
\begin{equation*}
  (khb \lor e)\inv{(b \lor e)} \in H \lor K. \qedhere 
\end{equation*}
\end{proof}

\end{proposition}

This conclusive result partially answers the Open Problem 6.28 presented in \cite{artalgebricallycoherent}. In fact, the category $\lab$ is algebraically coherent, as just shown, however it is not action accessible: the example provided in Proposition \ref{non accessibile} exclusively involves lattice-ordered groups whose group operations are commutative.

\section*{Acknowledgement}

The author would like to thank the referee for the useful comments and interesting suggestions, which have led to the present improved version of the paper. Gratitude is also extended to Andrea Montoli, the author's PhD supervisor, for his valuable insights and guidance, which were fundamental in the development of this work.

\vskip 1cm
\bibliography{Biblio}
\bibliographystyle{plain}
\end{document}